\documentclass[a4paper,12pt]{article}
\pdfoutput=1
\usepackage{etex}

\usepackage[inline]{asymptote}

\usepackage{amssymb}
\usepackage{amsmath}
\usepackage{amsthm}
\usepackage{bbm}

\usepackage{lineno}
\usepackage[cm]{fullpage}
\usepackage{graphicx}
\usepackage{wrapfig}
\usepackage{subfig}
\usepackage{float}
\usepackage{booktabs}
\usepackage{caption}

\allowdisplaybreaks

\makeatletter
\newsavebox{\sfe@box}
\newenvironment{subfloatenv}[1]{%
\def\sfe@caption{#1}%
\setbox\sfe@box\hbox\bgroup\color@setgroup}%
{\color@endgroup\egroup\subfloat[\sfe@caption]%
{\usebox{\sfe@box}}}
\makeatother

\theoremstyle{definition}

\newcommand{\norm}[1]{\ensuremath{\lVert#1\rVert}}

\newcommand{\reverse}[1]{\ensuremath{\widetilde{#1}}}

\newcommand{\R}[1]{\ensuremath{\mathbb{R}^{#1}}}

\newcommand{\E}[1]{\ensuremath{\mathbb{E}^{#1}}}
\newcommand{\M}[1]{\ensuremath{\mathbb{M}^{#1}}}

\newcommand{\El}[1]{\ensuremath{{E}_{#1}}}

\renewcommand{\S}[1]{\ensuremath{\mathbb{S}^{#1}}}

\newcommand{\tb}[1]{\ensuremath{\textbf{#1}}}
\newcommand{\mb}[1]{\ensuremath{\boldsymbol{#1}}}

\newcommand{\e}{\tb{e}} 
\newcommand{\I}{\tb{I}} 
\newcommand{\Id}{\ensuremath{I}} 
\newcommand{\J}{\ensuremath{J}} 

\newcommand{\p}{{\ensuremath{\!\!+}}}
\newcommand{\n}{{\ensuremath{\!\!-}}}

\title{Clifford algebra and the projective model of Elliptic spaces}
\author{Andrey Sokolov}

\begin{document}
\maketitle

\begin{abstract}
I apply the algebraic framework developed in \cite{gunn2011geometry} to study geometry of elliptic spaces in 1, 2, and 3 dimensions.
The background material on projectivised Clifford algebras and their application to Cayley-Klein geometries
is described in \cite{sokolov2013clifford}. 
The use of Clifford algebra largely obviates the need for spherical trigonometry as
elementary geometric transformations such as projections, rejections, reflections, and rotations can be
accomplished with geometric multiplication.
Furthermore, the same transformations can be used in 3-dimensional elliptic space where effective use of spherical trigonometry is problematic.
I give explicit construction of Clifford parallels and discuss their properties in detail.
Clifford translations are represented in a uniform fashion by geometric multiplication as well.
The emphasis of the exposition is on geometric structures and computation rather than proofs.  

\end{abstract}

\tableofcontents

\section{Elliptic line \El{1}}

In Elliptic space \El{1}, the norm of \(\tb{a}=d\e_0+a\e_1\) is given by \(\norm{\tb{a}}=\sqrt{d^2+a^2}\),
so any non-zero vector can be normalised.
A normalised vector can be written as 
\begin{equation}\label{parameterise El1}
\tb{a}=-\e_0\sin\alpha+\e_1\cos\alpha.
\end{equation}
It dually represents a point in \(E_1\) at \(x=\tan\alpha\) if \(\alpha\ne\tfrac{\pi}{2}+\pi k\) where \(k\in\mathbb{Z}\),
otherwise \(\tb{a}=\pm\e_0\).
Note that \(\e_1\tb{a}=\cos\alpha+\e_{01}\sin\alpha=e^{\alpha\e_{01}}\) and, therefore, \(\tb{a}\)
can be expressed via the exponential function as follows:
\begin{equation}
\tb{a}=\e_1e^{\alpha\e_{01}}.
\end{equation}
In the following, I identify vectors in the dual model space \R{2*} with the points  in \El{1} they dually represent. 
I will refer to \(\tb{a},\tb{b}\), etc as points in \El{1}.

\begin{figure}[h]
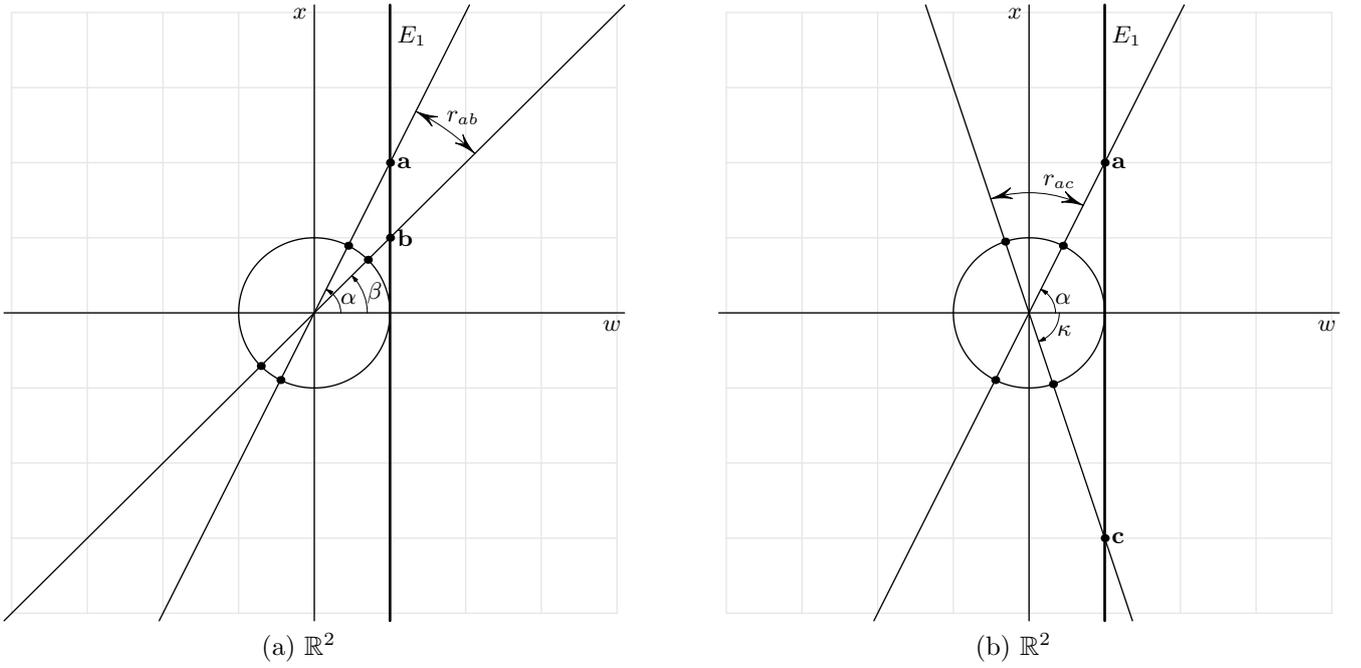

\begin{subfloatenv}{\R{2}}
\begin{asy}
import Figure2D;
Figure f = Figure(xaxis_name="$w$",yaxis_name="$x$");

real x(real alpha) { return sin(2*pi*alpha); }
real y(real alpha) { return cos(2*pi*alpha); }
path p = graph(x, y, 0, 1, operator ..);
draw(p);

f.line(e_1-e_0,draw_orientation=false,label="$E_1$",position=0.05,align=(1,0),pen=currentpen+1);

f.line(Line(0,-2,1),draw_orientation=false);
dot(Label("$\textbf{a}$"),(1,2));
f.arc1(O,Point(0,0,1),Point(0,-2,1),radius=0.35,label="$\alpha$",align=(0.5,0));
real alpha = atan(2);
dot((cos(alpha),sin(alpha)));
dot(-(cos(alpha),sin(alpha)));

f.line(Line(0,-1,1),draw_orientation=false);
dot(Label("$\textbf{b}$"),(1,1));
f.arc1(O,Point(0,0,1),Point(0,-1,1),radius=0.7,label="$\beta$",align=(0.5,0));
real beta = atan(1);
dot((cos(beta),sin(beta)));
dot(-(cos(beta),sin(beta)));

f.arc(O,Point(0,-1,1),Point(0,-2,1),radius=3,label="$r_{ab}$");

\end{asy}
\end{subfloatenv}\hfill%
\begin{subfloatenv}{\R{2}}
\begin{asy}
import Figure2D;
Figure f = Figure(xaxis_name="$w$",yaxis_name="$x$");

real x(real alpha) { return sin(2*pi*alpha); }
real y(real alpha) { return cos(2*pi*alpha); }
path p = graph(x, y, 0, 1, operator ..);
draw(p);

f.line(e_1-e_0,draw_orientation=false,label="$E_1$",position=0.05,align=(1,0),pen=currentpen+1);

f.line(Line(0,-2,1),draw_orientation=false);
dot(Label("$\textbf{a}$"),(1,2));
f.arc1(O,Point(0,0,1),Point(0,-2,1),radius=0.35,label="$\alpha$",align=(0.5,0));
real alpha = atan(2);
dot((cos(alpha),sin(alpha)));
dot(-(cos(alpha),sin(alpha)));

f.line(Line(0,3,1),draw_orientation=false);
dot(Label("$\textbf{c}$"),(1,-3));
f.arc1(O,Point(0,0,1),Point(0,3,1),direction=CW,radius=0.4,label="$\kappa$",align=(0.5,0));
real kappa = atan(-3);
dot((cos(kappa),sin(kappa)));
dot(-(cos(kappa),sin(kappa)));

f.arc(O,Point(0,-2,1),Point(0,-3,-1),radius=1.6,label="$r_{ac}$",align=(0.75,0.75));

\end{asy}
\end{subfloatenv}
\caption{The distance between points in \El{1}.}
\label{distance in El1}
\end{figure}

The distance \(r\in[0,\tfrac{\pi}{2}]\) between two normalised points \(\tb{a}\) and \(\tb{b}\)
is defined by
\begin{equation}
\sin r=|\tb{a}\vee\tb{b}|.
\end{equation}
In \El{1}, \(|\tb{a}\cdot\tb{b}|^2+|\tb{a}\vee\tb{b}|^2=1\) for the normalised points and, therefore,
the distance also satisfies \(\cos r=|\tb{a}\cdot\tb{b}|\).
For the normalised points 
\(\tb{a}=-\e_0\sin\alpha+\e_1\cos\alpha\) and \(\tb{b}=-\e_0\sin\beta+\e_1\cos\beta\),
I get \(\tb{a}\cdot\tb{b}=\sin\alpha\sin\beta+\cos\alpha\cos\beta=\cos(\alpha-\beta)\)
and, therefore,
\begin{equation}
r=
\left\{
\begin{aligned}
&|\alpha-\beta|, \textrm{ if } |\alpha-\beta|\le\tfrac{\pi}{2},\\
&\pi-|\alpha-\beta|, \textrm{ if } \tfrac{\pi}{2}<|\alpha-\beta|\le\pi.
\end{aligned}
\right.
\end{equation}
The distance in \El{1} has a maximum value of \(\tfrac{\pi}{2}\).
For instance, the distance from the origin \(\e_1\) to \(\e_0\) equals \(\tfrac{\pi}{2}\),
so \(\e_0\) is not at infinity in \El{1}.
In fact, \(\e_0\) is not different from any other point in \El{1}. 

The points \(\tb{a}=\tfrac{1}{\sqrt{5}}(-2\e_0+\e_1)\) and  \(\tb{b}=\tfrac{1}{\sqrt{2}}(-\e_0+\e_1)\)
shown in Figure~\ref{distance in El1}(a) are parametrised by
\(\alpha=\arctan{2}\), \(\beta=\tfrac{\pi}{4}\) via (\ref{parameterise El1}), which
gives \(r_{ab}=\arctan{2}-\tfrac{\pi}{4}\) for the distance between them.
Similarly, for the points \(\tb{a}=\tfrac{1}{\sqrt{5}}(-2\e_0+\e_1)\) and \(\tb{c}=\tfrac{1}{\sqrt{10}}(3\e_0+\e_1)\)
shown in Figure~\ref{distance in El1}(b),
I get \(r_{ac}=\pi-(\arctan{2}+\arctan{3})\) for the distance between them
(\(\tb{c}\) is parametrised by \(\kappa=-\arctan{3}\)).
In both cases, the distance between the points 
is equal to the angle (in Euclidean sense in \R{2}) between the linear subspaces representing the points.
Due to this property of the distance measure in Elliptic space,
it is convenient to visualise  \El{1} as a unit circle centered on the origin 
as shown in Figure~\ref{distance in El1},
instead of a straight line at \(w=1\).
The circle is defined by 
\begin{equation}
w^2+x^2=1
\end{equation}
and the correspondence between the linear subspaces of \R{2} and the points on the circle is 
determined by the intersection of the subspaces with the circle.
A linear subspace intersects the circle in two points, which by construction are identified and treated as a single point in \El{1}.
In this representation of \El{1}, it is almost equivalent to a 1-dimensional version of spherical geometry denoted by \S{1}.
The difference between \El{1} and \S{1} is that the antipodal points on the circle,
which are treated as separate points in \S{1}, are identical in \El{1}.

The polar point \(\tb{a}\e_{01}\) of a point \(\tb{a}\) is a point in \El{1} 
which is at a distance of \(\tfrac{\pi}{2}\) from \(\tb{a}\).
This can be readily seen from 
\(\tb{a}\e_{01}
=\norm{a}\e_1e^{\alpha\e_{01}}e^{\tfrac{\pi}{2}\e_{01}}
=\norm{a}\e_1 e^{(\alpha+\tfrac{\pi}{2})\e_{01}}\),
where I used \(\e_{01}=e^{\tfrac{\pi}{2}\e_{01}}\).
For instance, \(\e_0\) is the polar point of \(-\e_1\), and \(\e_1\) is the polar point of \(\e_0\).

Spinors in \El{1} are defined in the standard way as the product of an even number of normalised proper points.
Note that every point in \El{1} is proper, including \(\e_0\).
Any spinor can be written as \(S=e^{\alpha\e_{01}}\) for some \(\alpha\in\R{}\), e.g.\ \((-\e_1)\e_1=-1=e^{\pi\e_{01}}\).
Spinors in \El{1} can be identified with normalised complex numbers.

The translation of a point \(\tb{a}\) by \(\lambda\) is given by \(T\tb{a}T^{-1}\), 
where \(T=e^{-\tfrac{1}{2}\lambda\e_{01}}\).
If \(\tb{a}\) is parametrised with \(\alpha\), then \(\tb{a}'=T\tb{a}T^{-1}\) can be parametrised with 
\(\alpha'=\alpha+\lambda\), which follows from 
\(\tb{a}'=T\tb{a}T^{-1}=\tb{a}T^{-2}
=\e_1  e^{\alpha\e_{01}}e^{\lambda\e_{01}}
=\e_1e^{(\alpha+\lambda)\e_{01}}\) where I used \(\tb{a}=\e_1e^{\alpha\e_{01}}\).
Note that \(\tb{a}=\e_1e^{\alpha\e_{01}}\) can be interpreted as the translation of \(\e_1\), the origin of \El{1}, by \(\alpha\).
So, a translation in \El{1} is equivalent to a Euclidean rotation around the origin of \R{2}.
The space \El{1} is closed and periodic with the period equal to \(\pi\);
translating any point by \(\pi\) yields the original point (albeit with the opposite orientation).

The top-down reflection of \(\tb{a}\) in \(\tb{b}\) is given by \(-\tb{b}\tb{a}\tb{b}^{-1}\).
If \(\tb{a}\) is parametrised with \(\alpha\) and \(\tb{b}\) with \(\beta\), then
\(\tb{a}'=-\tb{b}\tb{a}\tb{b}^{-1}\) can be parametrised with \(\alpha'=\beta-(\alpha-\beta)\), which follows from
\(\tb{a}'=-\tb{b}\tb{a}\tb{b}^{-1}=-\e_1e^{\beta\e_{01}} \e_1e^{\alpha\e_{01}} e^{-\beta\e_{01}}\e_1^{-1}
=-\e_1e^{\beta\e_{01}} e^{-\alpha\e_{01}} e^{\beta\e_{01}}
=-\e_1e^{(\beta-\alpha+\beta)\e_{01}}\).
Note that the orientation of \(\tb{a}'\) is the opposite of \(\tb{a}\).

The projection of \(\tb{a}\) on \(\tb{b}\) is given by the familiar expression
\((\tb{a}\cdot\tb{b})\tb{b}^{-1}\),
which yields \(\tb{b}\cos(\alpha-\beta)\) if \(\tb{a}\) and \(\tb{b}\) are normalised
(\(\tb{a}\) is parametrised by \(\alpha\), and \(\tb{b}\)  by \(\beta\)).
Indeed, 
\((\tb{a}\cdot\tb{b})\tb{b}^{-1}=
\tfrac{1}{2}(\tb{a}\tb{b}+\tb{b}\tb{a})\tb{b}
=\tfrac{1}{2}\tb{b}(\tb{b}\tb{a}+\tb{a}\tb{b})
=\tb{b}\tfrac{1}{2}(e^{(\alpha-\beta)\e_{01}}+e^{(\beta-\alpha)\e_{01}})
=\tb{b}\cos(\alpha-\beta)
\).
So, if \(\tb{a}\) is at a distance of \(\tfrac{\pi}{2}\) from \(\tb{b}\),
the projection of \(\tb{a}\) on \(\tb{b}\) is zero.
Unlike \(\e_0\) in Euclidean space, \(\e_0\) in \El{1} is invertible and it is possible to project on \(\e_0\) as on any other point.
The rejection of \(\tb{a}\) by \(\tb{b}\) is given by \((\tb{a}\wedge\tb{b})\tb{b}^{-1}\),
which for the normalised points yields \(\tb{b}\e_{01}\sin(\alpha-\beta)\).
So, the rejection by \(\tb{b}\) coincides with the polar point of \(\tb{b}\) weighted by \(\sin(\alpha-\beta)\).
Moreover, the rejection is equivalent to projection on the polar point.
The expression for scaling, which I used in degenerate spaces, is not applicable in a non-degenerate space.


\section{Elliptic plane \El{2}}
The norm of  a line \(\tb{a}=d\e_0+a\e_1+b\e_2\) and a point \(\tb{P}=w\e_{12}+x\e_{20}+y\e_{01}\) 
 is given by \(\norm{\tb{a}}=\sqrt{d^2+a^2+b^2}\) and \(\norm{\tb{P}}=\sqrt{w^2+x^2+y^2}\).
So, any line or point in \El{2} can be normalised.

\begin{figure}[t!]
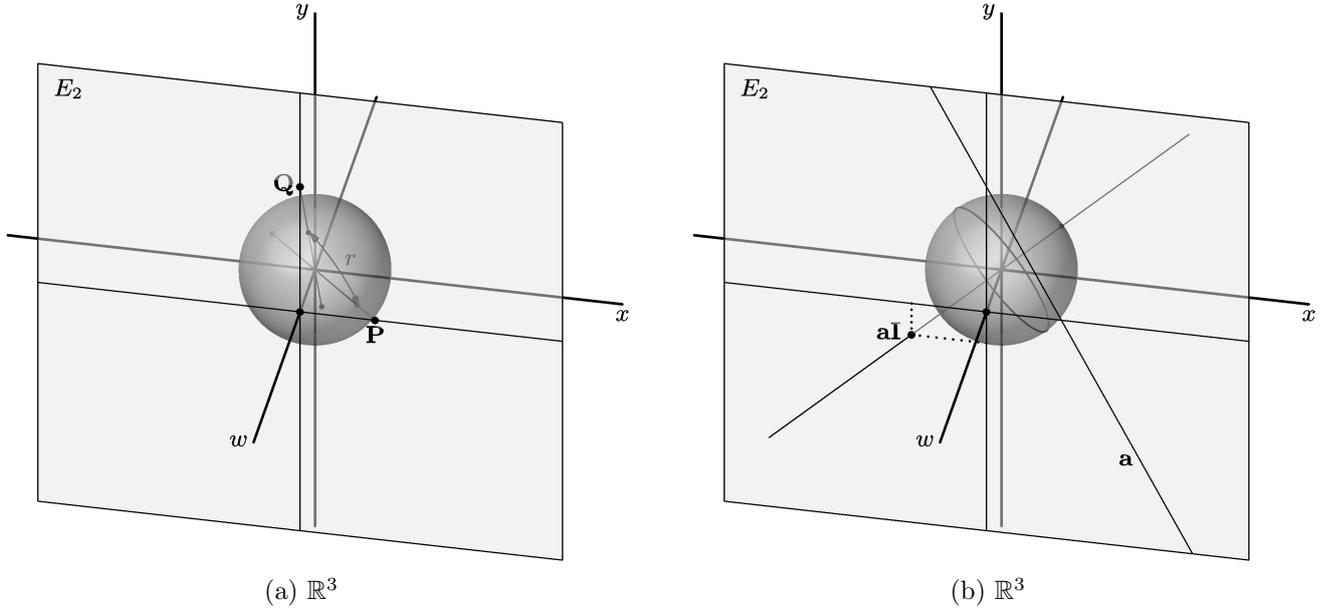

\begin{subfloatenv}{\R{3}}
\begin{asy}
import Drawing3D;
DrawingR3 drawing = DrawingR3(4.0,0.1, camera=(10,2,7));	
drawing.target_axes();

path3 p = plane_path3((1,0,0),(1,0,0), v=(0,1,0));
drawing.Drawing3D.plane(p);
dot((1,0,0));
draw((1,-3.5,0)--(1,3.5,0));
draw((1,0,-3.5)--(1,0,3.5));

label("$E_2$", (1.5,-3,3.5));

import solids;
real r = 1;
currentlight=light(white, viewport=false,(5,-5,10));
draw(surface(sphere((0,0,0),r,100)),lightgray+opacity(0.4));

triple P=(1,1,0);
triple Q=(1,0,2);

label("$\textbf{P}$",P,(0,0,-1));
label("$\textbf{Q}$",Q,(0,-1,0));
pen dot_pen = currentpen+2;
dot(P); dot(normalise(P), dot_pen); dot(normalise(-P), dot_pen);
dot(Q); dot(normalise(Q), dot_pen); dot(normalise(-Q), dot_pen);
draw(normalise(-P)--P);
draw(normalise(-Q)--Q);

path3 PQ = arc((0,0,0), 1.05*normalise(P), 1.05*normalise(Q));
draw(Label("$r$",(1,1.5,2)),PQ, Arrows3(size=5));
\end{asy}
\end{subfloatenv}\hspace{-12pt}%
\begin{subfloatenv}{\R{3}}
\begin{asy}
import Drawing3D;
DrawingR3 drawing = DrawingR3(4.0,0.1, camera=(10,2,7));	
drawing.target_axes();

path3 p = plane_path3((1,0,0),(1,0,0), v=(0,1,0));
drawing.Drawing3D.plane(p);
dot((1,0,0));
draw((1,-3.5,0)--(1,3.5,0));
draw((1,0,-3.5)--(1,0,3.5));

path3 Q = plane_path3((1,-1,-1/2), (1,1,0), theta=-13, scale=7.8);

path3 l = plane_intersection(Q,p);
draw(Label("$\textbf{a}$", 0.2,W), l);

label("$E_2$", (1.5,-3,3.5));

import solids;
real r = 1;
currentlight=light(white, viewport=false,(5,-5,10));
draw(surface(sphere((0,0,0),r,100)),lightgray+opacity(0.4));
path3 c = Circle((0,0,0), 1, (1,-1,-1/2));
draw(c,black+linewidth(0.25mm));

path3 l = line_path3((1,-1,-1/2)/4, direction=(1,-1,-1/2), extent_from_center=3.5);
draw(l);
label("$\textbf{a}\textbf{I}$",(1,-1,-1/2),(1/2,-1-1/2,1/2));
dot((1,-1,-1/2));
pen dot_pen = currentpen+2;
dot(normalise((1,-1,-1/2)), dot_pen);
dot(normalise(-(1,-1,-1/2)), dot_pen);

draw((1,-1,0)--(1,-1,-1/2), Dotted);
draw((1,0,-1/2)--(1,-1,-1/2), Dotted);

\end{asy}
\end{subfloatenv}
\caption{Visualising \El{2}}
\label{visualising elliptic2}
\end{figure}

The distance \(r\in[0,\tfrac{\pi}{2}]\) between two normalised points \(\tb{P}\) and \(\tb{Q}\) is defined by
\begin{equation}\label{distance in El2}
\sin r=\norm{\tb{P}\vee\tb{Q}}.
\end{equation}
Since \(|\tb{P}\cdot\tb{Q}|^2+\norm{\tb{P}\vee\tb{Q}}^2=1\) for normalised points in  \El{2}, the distance \(r\) also satisfies
\(\cos r=|\tb{P}\cdot\tb{Q}|\).
If the point \(\tb{P}\) is expressed via the standard coordinates \(x\) and \(y\), 
i.e.\ \(\tb{P}=\e_{12}+x\e_{20}+y\e_{01}\), then \(\norm{\tb{P}}=\sqrt{1+x^2+y^2}\) and
\begin{equation}
\sin{r}=\frac{\sqrt{x^2+y^2}}{\sqrt{1 + x^2+y^2}}, \quad\quad
\cos{r}=\frac{1}{\sqrt{1 + x^2+y^2}}.
\end{equation}
So, those points that are at an infinite distance from the origin in Euclidean space \E{2} are
at the distance of \(\tfrac{\pi}{2}\) from the origin in Elliptic space \El{2}, which is the maximum distance in \El{2}.
%
The distance between arbitrary points \(\tb{P}\) and \(\tb{Q}\) 
can be understood as the Euclidean angle between the 1-dimensional linear subspaces in \R{3} 
corresponding to \(\J(\tb{P})\) and \(\J(\tb{Q})\), where \R{3} is viewed as a 3-dimensional Euclidean space.
See Figure~\ref{visualising elliptic2}(a), where \(\tb{P}=\e_{12}+\e_{20}\) and \(\tb{Q}=\e_{12}+2\e_{01}\), 
which upon the normalisation
yields \(r=\arccos(\tfrac{1}{\sqrt{2}}\tfrac{1}{\sqrt{5}})\) for the distance between the points.

The angle \(\alpha\in[0,\pi]\) between two normalised lines \(\tb{a}\) and \(\tb{b}\),
or more precisely the angle between the orientation vectors of the lines, is defined by
\begin{equation}
\cos\alpha=\tb{a}\cdot\tb{b}.
\label{angle in E2}
\end{equation}
It corresponds to the Euclidean angle between \(\tb{a}\) and \(\tb{b}\) only if both lines pass through the origin.
In general, \(\alpha\) is equal to the Euclidean angle between the 2-dimensional linear subspaces in \R{3} 
corresponding to \(\J(\tb{a})\) and \(\J(\tb{b})\),
or more precisely it is equal to the Euclidean angle between the vectors \(I(\tb{a})\) and \(I(\tb{b})\),
which are perpendicular (in the Euclidean sense of \R{3}) to the planes \(\J(\tb{a})\) and \(\J(\tb{b})\).
Since both the distance and angular measures are elliptic, there is a close relationship between them.
For instance, if \(0\le\alpha\le\tfrac{\pi}{2}\), then \(\alpha\) is equal to the distance between the polar points
of the lines  \(\tb{a}\) and \(\tb{b}\), i.e. \(\cos \alpha=|(\tb{a}\I)\cdot(\tb{b}\I)|\).

Due to the above mentioned properties of the distance and angular measures,
it might be more convenient to visualise \El{2} as a unit sphere centered on the origin in \R{3} rather than 
a plane at \(w=1\).
The sphere is defined by
\begin{equation}
w^2+x^2+y^2=1
\end{equation}
and its antipodal points are identified.
This is a natural extension of the representation of \El{1} as a unit circle.

When \El{2} is viewed as the unit sphere, each line in \El{2} can be identified with a great circle of the sphere.
For instance, the line \(\tb{a}=\tb{P}\vee\tb{Q}=-2\e_0+2\e_1+\e_2\) shown in Figure~\ref{visualising elliptic2}(b)
is represented by a 2-dimensional linear subspace of \R{3}, 
which intersects the unit sphere along a great circle.
If the great circle is thought of as the equator, 
then the polar point \(\tb{a}\I=-2\e_{12}+2\e_{20}+\e_{01}=-2(\e_{12}-\e_{20}-\tfrac{1}{2}\e_{01})\) of the line \(\tb{a}\) 
is located at the poles of the sphere with respect to the equator.
The polar point \(\tb{a}\I\) is at the distance of \(\tfrac{\pi}{2}\) from any point on the line \(\tb{a}\).
Hence, the commutator \(\tb{P}\times\tb{Q}\), which coincides with the polar point of the line \(\tb{P}\vee\tb{Q}\),
is located at a finite distance from the origin.
Since \(\I^{-1}=-\I\) in \El{2}, the definition of the duality transformation implies that  \(\J(\tb{a})=-\Id(\tb{a}\I)\),
where \(\J\) and \(\Id\) are the duality and identity transformations.
So, a line \(\tb{a}\) and its polar point \(\tb{a}\I\) are related to each other in \El{2} in essentially the same
way as a line and a point dual to each other in the sense of projective duality.

\begin{figure}[t!]
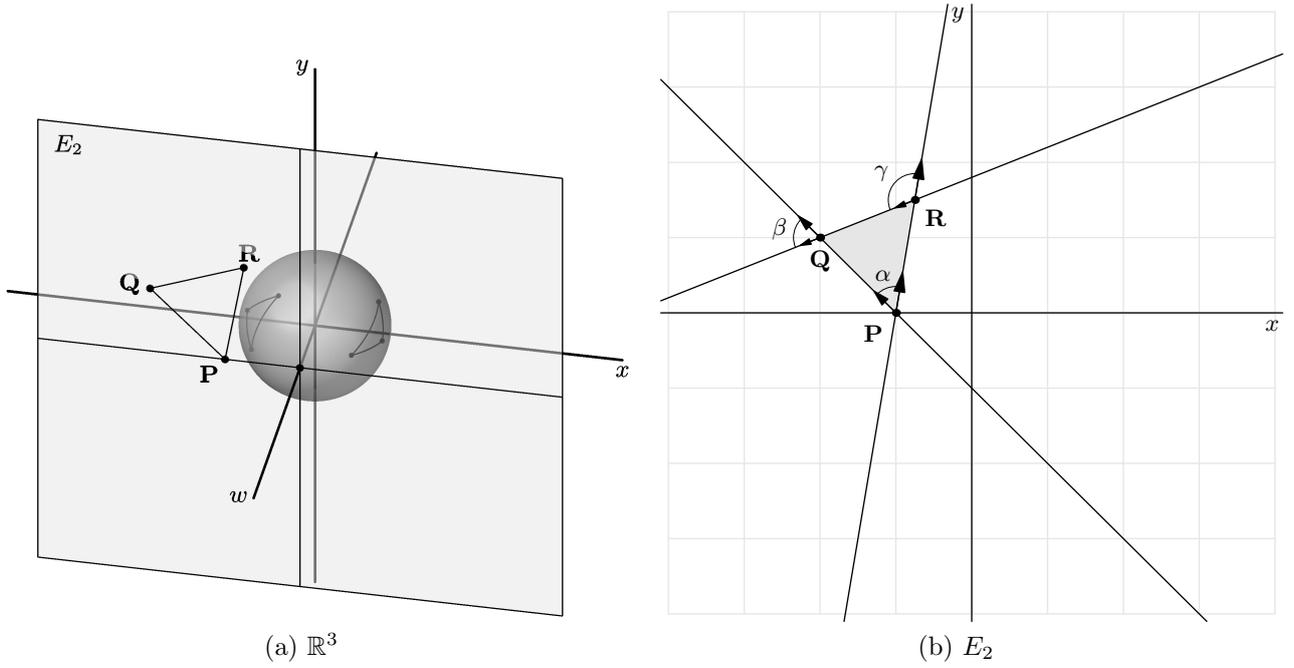

\begin{subfloatenv}{\R{3}}
\begin{asy}
import Drawing3D;
DrawingR3 drawing = DrawingR3(4.0,0.1, camera=(10,2,7));	
drawing.target_axes();

path3 p = plane_path3((1,0,0),(1,0,0), v=(0,1,0));
drawing.Drawing3D.plane(p);
dot((1,0,0));
draw((1,-3.5,0)--(1,3.5,0));
draw((1,0,-3.5)--(1,0,3.5));

path3 Q = plane_path3((1,-1,-1/2), (1,1,0), theta=-13, scale=7.8);
label("$E_2$", (1.5,-3,3.5));

import solids;
real r = 1;
currentlight=light(white, viewport=false,(5,-5,10));
draw(surface(sphere((0,0,0),r,100)),lightgray+opacity(0.4));

import Projective2D;
metric=Metric(Elliptic);
MV P = Point(1,-1,0); P/=norm(P);
MV Q =  Point(1,-2,1); Q/=norm(Q);
MV R = Point(1,-3/4,3/2); R/=norm(R);

triple toxy(MV P) {return (1,topair(P).x,topair(P).y); }

dot(toxy(P)); label("$\textbf{P}$",toxy(P),(1/2,-1,-1));
dot(toxy(Q)); label("$\textbf{Q}$",toxy(Q),(1/2,-1-1/2,1/2));
dot(toxy(R)); label("$\textbf{R}$",toxy(R),(1/2,1/2,1.75));

draw(toxy(P)--toxy(Q));
draw(toxy(P)--toxy(R));
draw(toxy(R)--toxy(Q));

triple to3(MV point) {return (point.w, point.x, point.y);};
pen dot_pen = currentpen+2;
dot(to3(P), dot_pen); dot(-to3(P), dot_pen);
dot(to3(Q), dot_pen); dot(-to3(Q), dot_pen);
dot(to3(R), dot_pen); dot(-to3(R), dot_pen);

triple F(real t, MV P1, MV P2) { MV R=-join(P1,P2)*I; R/=norm(R); MV S = cos(1/2*t)-R*sin(1/2*t); return to3(S*P1/S); };
real span(MV P1, MV P2) { real c=dot(P1,P2).s; if(c<0) return pi-acos(c); else return acos(c); };

triple f(real t) {return F(t,P,Q);};
real s = span(P,Q);
real x(real t) { return f(t).x; };
real y(real t) { return f(t).y; };
real z(real t) { return f(t).z; };
path3 a1 = graph(x, y, z, 0, s, operator ..);
draw(a1);
path3 a1 = graph(x, y, z, pi, pi+s, operator ..);
draw(a1);

triple f(real t) {return F(t,Q,R);};
real s = span(Q,R);
real x(real t) { return f(t).x; };
real y(real t) { return f(t).y; };
real z(real t) { return f(t).z; };
path3 b1 = graph(x, y, z, 0, s, operator ..);
draw(b1);
path3 b1 = graph(x, y, z, pi, pi+s, operator ..);
draw(b1);

triple f(real t) {return F(t,R,P);};
real s = span(R,P);
real x(real t) { return f(t).x; };
real y(real t) { return f(t).y; };
real z(real t) { return f(t).z; };
path3 c1 = graph(x, y, z, 0, s, operator ..);
draw(c1);
path3 c1 = graph(x, y, z, pi, pi+s, operator ..);
draw(c1);

\end{asy}
\end{subfloatenv}\hspace{-12pt}%
\begin{subfloatenv}{\El{2}}
\begin{asy}
import Figure2D;
Figure f = Figure();
metric=Metric(Elliptic);
MV P = Point(1,-1,0); P/=norm(P);
MV Q = Point(1,-2,1); Q/=norm(Q);
MV R = Point(1,-3/4,3/2); R/=norm(R);

fill(topair(P/P.w)--topair(Q/Q.w)--topair(R/R.w)--cycle,lightgray);
f.point(P,"$\textbf{P}$",align=(-1,-1),draw_orientation=false);
f.point(Q,"$\textbf{Q}$",align=(0,-1),draw_orientation=false);
f.point(R,"$\textbf{R}$",align=(0.8,-0.8),draw_orientation=false);

MV rq = join(R,Q);
MV pr = join(P,R);
MV pq = join(P,Q);

f.line(pq,O=topair(P),draw_orientation=false,draw_bottom_up_orientation=true);
f.line(pr,O=topair(P),draw_orientation=false,draw_bottom_up_orientation=true);
f.arc2(P,wedge(e_0,pr),wedge(e_0,pq),"$\alpha$");

f.line(rq,pen=invisible,O=topair(Q),draw_orientation=false,draw_bottom_up_orientation=true);
f.line(pq,pen=invisible,O=topair(Q),draw_orientation=false,draw_bottom_up_orientation=true);
f.arc2(Q,wedge(e_0,pq),wedge(e_0,rq),"$\beta$");
 
f.line(rq,O=topair(R),draw_orientation=false,draw_bottom_up_orientation=true);
f.line(pr,pen=invisible,O=topair(R),draw_orientation=false,draw_bottom_up_orientation=true);
f.arc2(R,wedge(e_0,pr),wedge(e_0,rq), "$\gamma$");

rq/=norm(rq);
pr/=norm(pr);
pq/=norm(pq);

real alpha = acos(dot(pr,pq).s);
real beta = acos(dot(pq,rq).s);
real gamma = acos(dot(pr,rq).s);
alpha+beta-gamma;

\end{asy}
\end{subfloatenv}
\caption{Triangles in \El{2} (1)}
\label{triangles in El2}
\end{figure}%
\begin{figure}[t!]
\begin{subfloatenv}{\R{3}}
\begin{asy}
import Drawing3D;
DrawingR3 drawing = DrawingR3(4.0,0.1, camera=(10,2,7));	
drawing.target_axes();

path3 p = plane_path3((1,0,0),(1,0,0), v=(0,1,0));
drawing.Drawing3D.plane(p);
dot((1,0,0));
real bb = 3.5;
draw((1,-bb,0)--(1,bb,0));
draw((1,0,-bb)--(1,0,bb));

path3 Q = plane_path3((1,-1,-1/2), (1,1,0), theta=-13, scale=7.8);
label("$E_2$", (1.5,-bb+1,bb));

import solids;
real r = 1;
currentlight=light(white, viewport=false,(5,-5,10));
draw(surface(sphere((0,0,0),r,100)),lightgray+opacity(0.4));

import Projective2D;
metric=Metric(Elliptic);
MV P = Point(1,3,-1); P/=norm(P);
MV Q = Point(1,-2,1); Q/=norm(Q);
MV R = Point(1,-3/4,3/2); R/=norm(R);

triple toxy(MV P) {return (1,topair(P).x,topair(P).y); }

dot(toxy(P)); label("$\textbf{P}$",toxy(P),(1/2,-1,0));
dot(toxy(Q)); label("$\textbf{Q}$",toxy(Q),(1/2,-0,1.5));
dot(toxy(R)); label("$\textbf{R}$",toxy(R),(1/2,1/2,1.75));

triple to3(MV point) {return (point.w, point.x, point.y);};
pen dot_pen = currentpen+2;
dot(to3(P), dot_pen); dot(-to3(P), dot_pen);
dot(to3(Q), dot_pen); dot(-to3(Q), dot_pen);
dot(to3(R), dot_pen); dot(-to3(R), dot_pen);

MV PQ1 = wedge(join(P,Q),Line(bb,1,0));
MV PQ2 = wedge(join(P,Q),Line(bb,-1,0));

MV PR1 = wedge(join(P,R),Line(bb,1,0));
MV PR2 = wedge(join(P,R),Line(bb,-1,0));

draw(toxy(Q)--toxy(PQ1)); draw(toxy(P)--toxy(PQ2));
draw(toxy(R)--toxy(PR1)); draw(toxy(P)--toxy(PR2));
draw(toxy(R)--toxy(Q));

triple F(real t, MV P1, MV P2) { MV R=-join(P1,P2)*I; R/=norm(R); MV S = cos(1/2*t)-R*sin(1/2*t); return to3(S*P1/S); };
real span(MV P1, MV P2) { real c=dot(P1,P2).s; if(c<0) return pi-acos(c); else return acos(c); };

triple f(real t) {return -F(-t,P,Q);};
real s = span(P,Q);
real x(real t) { return f(t).x; };
real y(real t) { return f(t).y; };
real z(real t) { return f(t).z; };
path3 a1 = graph(x, y, z, 0, s, operator ..);
draw(a1);
path3 a1 = graph(x, y, z, pi, pi+s, operator ..);
draw(a1);

triple f(real t) {return F(t,Q,R);};
real s = span(Q,R);
real x(real t) { return f(t).x; };
real y(real t) { return f(t).y; };
real z(real t) { return f(t).z; };
path3 b1 = graph(x, y, z, 0, s, operator ..);
draw(b1);
path3 b1 = graph(x, y, z, pi, pi+s, operator ..);
draw(b1);

triple f(real t) {return -F(-t,R,P);};
real s = span(R,P);
real x(real t) { return f(t).x; };
real y(real t) { return f(t).y; };
real z(real t) { return f(t).z; };
path3 c1 = graph(x, y, z, 0, s, operator ..);
draw(c1);
path3 c1 = graph(x, y, z, pi, pi+s, operator ..);
draw(c1);

\end{asy}
\end{subfloatenv}\hspace{-12pt}%
\begin{subfloatenv}{\El{2}}
\begin{asy}
import Figure2D;
Figure f = Figure();
metric=Metric(Elliptic);
MV P = Point(1,3,-1); P/=norm(P);
MV Q = Point(1,-2,1); Q/=norm(Q);
MV R = Point(1,-3/4,3/2); R/=norm(R);

real bb = f.full_extent_from_origin;
MV PQ1 = wedge(join(P,Q),Line(bb,1,0));
MV PQ2 = wedge(join(P,Q),Line(bb,-1,0));

MV PR1 = wedge(join(P,R),Line(bb,1,0));
MV PR2 = wedge(join(P,R),Line(bb,-1,0));

fill(topair(P/P.w)--topair(PQ2/PQ2.w)--topair(PR2/PR2.w)--cycle,lightgray);
fill(topair(Q/Q.w)--topair(R/R.w)--topair(PR1/PR1.w)--topair(PQ1/PQ1.w)--cycle,lightgray);

f.point(P,"$\textbf{P}$",align=(-1,-1),draw_orientation=false);
f.point(Q,"$\textbf{Q}$",align=(0,1),draw_orientation=false);
f.point(R,"$\textbf{R}$",align=(0,3/4),draw_orientation=false);

MV rq = join(R,Q);
MV pr = -join(P,R);
MV pq = -join(P,Q);

f.line(pq,O=topair(P),draw_orientation=false,draw_bottom_up_orientation=true);
f.line(pr,O=topair(P),draw_orientation=false,draw_bottom_up_orientation=true);
f.arc2(P,wedge(e_0,pr),wedge(e_0,pq),"$\alpha$",radius=1,align=(-1,-2));

f.line(rq,pen=invisible,O=topair(Q),draw_orientation=false,draw_bottom_up_orientation=true);
f.line(pq,pen=invisible,O=topair(Q),draw_orientation=false,draw_bottom_up_orientation=true);
f.arc2(Q,wedge(e_0,rq),wedge(e_0,pq),"$\beta$",radius=0.25);
 
f.line(rq,O=topair(R),draw_orientation=false,draw_bottom_up_orientation=true);
f.line(pr,pen=invisible,O=topair(R),draw_orientation=false,draw_bottom_up_orientation=true);
f.arc2(R,wedge(e_0,rq),wedge(e_0,pr),"$\gamma$",radius=0.25);

rq/=norm(rq);
pr/=norm(pr);
pq/=norm(pq);

real alpha = acos(dot(pr,pq).s);
real beta = acos(dot(pq,rq).s);
real gamma = acos(dot(rq,pr).s);
alpha+beta-gamma;

\end{asy}
\end{subfloatenv}
\caption{Triangles in \El{2} (2)}
\label{triangles in El2 unusual}
\end{figure}

The distance \(r\) between a normalised line \(\tb{a}\) and a normalised point \(\tb{P}\) satisfies
\(\sin r=|\tb{a}\vee\tb{P}|\) and alternatively \(\cos r=\norm{\tb{a}\cdot\tb{P}}\).
For instance, the distance from \(\tb{P}=\e_{12}-\tfrac{3}{5}\e_{20}+\tfrac{4}{5}\e_{01}\) to 
\(\tb{a}=-2\e_0+ 2\e_1 +\e_2\) is given by \(\sin r=\tfrac{2.4}{3\sqrt{2}}\),
where normalisation has been taken into account.
For a normalised line \(\tb{a}=d\e_0+a\e_1+b\e_2\), I get \(\sin r=|d|\) for the distance \(r\)
from the line to the origin.
For instance, the distance of \(\e_0\) from the origin equals \(\tfrac{\pi}{2}\).

Unlike the flat spaces \E{2} and \M{2},
there is no one-to-one correspondence between the area \({\cal S}\) of a triangle in \El{2}
defined by three normalised points \(\tb{P}\), \(\tb{Q}\), \(\tb{R}\) and the value of \(|\tb{P}\vee\tb{Q}\vee\tb{R}|\).
However, if one of the angles of the triangle is right, say, the angle at the point \(\tb{P}\) equals \(\tfrac{\pi}{2}\),
then the area \({\cal S}\) satisfies
\begin{equation}
\sin{\cal S}=\frac{|\tb{P}\vee\tb{Q}\vee\tb{R}|}{1+|\tb{Q}\cdot\tb{R}|},
\label{area in El2 right}
\end{equation}
where \(\tb{P}\), \(\tb{Q}\), and \(\tb{R}\) are assumed to be normalised and 
\((\tb{P}\vee\tb{Q})\cdot(\tb{P}\vee\tb{R})=0\).
In particular, if all three points are at the distance \(\tfrac{\pi}{2}\) from one another,
then \(\sin{\cal S}=1\) and \({\cal S}=\tfrac{\pi}{2}\), which is the maximum possible area of a triangle in \El{2}. 
The total area of \El{2} equals \(2\pi\).

An example of a triangle is shown in Figure~\ref{triangles in El2}, where
\(\tb{P}=\tfrac{1}{\sqrt{2}}(\e_{12}-\e_{20})\),
\(\tb{Q}=\tfrac{1}{\sqrt{6}}(\e_{12}-2\e_{20}+\e_{01})\), and
\(\tb{R}=\tfrac{4}{\sqrt{61}}(\e_{12}-\tfrac{3}{4}\e_{20}+\tfrac{3}{2}\e_{01})\).
The angles associated with the triangle (see Figure~\ref{triangles in El2}(b))
can be found by applying the definition~(\ref{angle in E2}), which gives
\(\alpha=\arccos(\frac{\tb{r}\cdot\tb{q}}{\norm{\tb{r}}\norm{\tb{q}}})\),
\(\beta=\arccos(\frac{\tb{r}\cdot\tb{p}}{\norm{\tb{r}}\norm{\tb{p}}})\),
\(\gamma=\arccos(\frac{\tb{q}\cdot\tb{p}}{\norm{\tb{q}}\norm{\tb{p}}})\),
where \(\tb{r}=\tb{P}\vee\tb{Q}\), \(\tb{q}=\tb{P}\vee\tb{R}\), \(\tb{p}=\tb{R}\vee\tb{Q}\).
Observe the correspondence between the angles as shown in Figure~\ref{triangles in El2}(b)
  and the bottom-up orientation of lines,
which constitute the sides of the triangle.
Then, the area \(\cal S\) can be computed with 
\begin{equation}
\cal S=\alpha+\beta-\gamma,
\end{equation}
which is equivalent to the spherical excess formula \(\cal S=\alpha+\beta+(\pi-\gamma)-\pi\).
The same area can be computed by breaking the triangle into two right triangles and
computing the area of each by applying~(\ref{area in El2 right}).

Another example is shown in Figure~\ref{triangles in El2 unusual}, where
\(\tb{Q}\) and \(\tb{R}\) are the same and 
\(\tb{P}=\tfrac{1}{\sqrt{11}}(\e_{12}+3\e_{20}-\e_{01})\).
Note the unusual appearance of this triangle in the plane \El{2}.
The lines by the sides of the triangle are defined by 
\(\tb{r}=\tb{Q}\vee\tb{P}\), \(\tb{q}=\tb{R}\vee\tb{P}\), \(\tb{p}=\tb{R}\vee\tb{Q}\),
which gives the bottom-up orientation shown in Figure~\ref{triangles in El2 unusual}(b).
Then, \(\cal S=\alpha+\beta-\gamma\) as before.

\begin{figure}[t!]
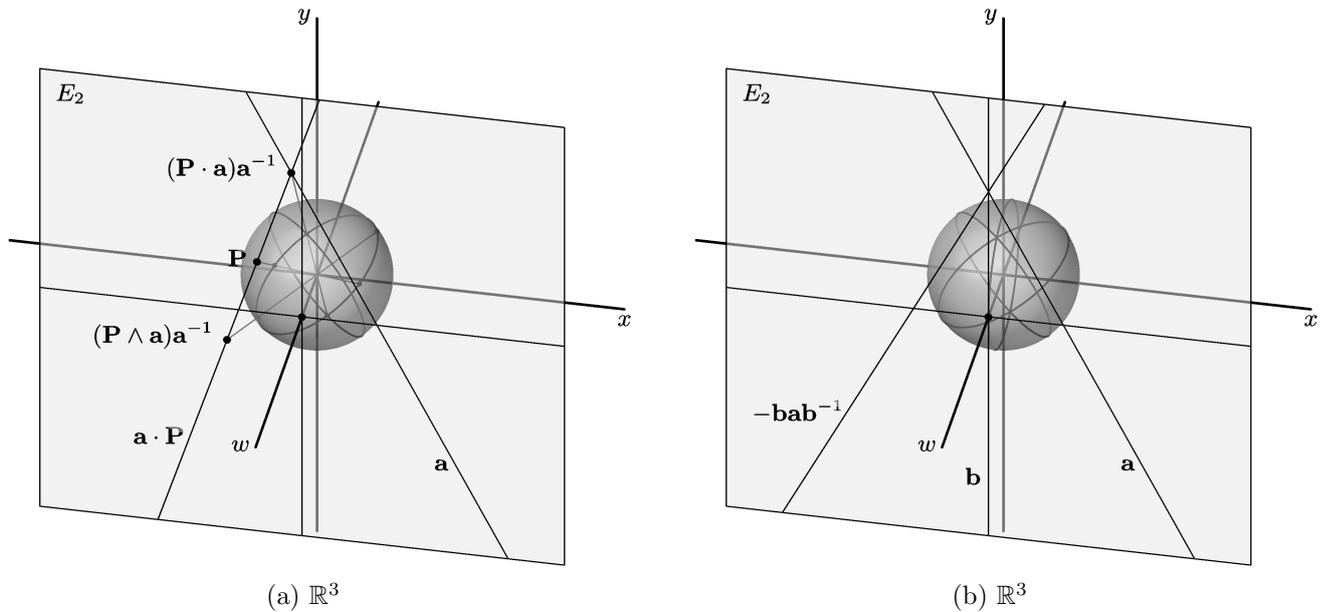

\begin{subfloatenv}{\R{3}}
\begin{asy}
import Drawing3D;
DrawingR3 drawing = DrawingR3(4.0,0.1, camera=(10,2,7));	
drawing.target_axes();

path3 p = plane_path3((1,0,0),(1,0,0), v=(0,1,0));
drawing.Drawing3D.plane(p);
dot((1,0,0));
draw((1,-3.5,0)--(1,3.5,0));
draw((1,0,-3.5)--(1,0,3.5));

path3 Q = plane_path3((1,-1,-1/2), (1,1,0), theta=-13, scale=7.8);

path3 l = plane_intersection(Q,p);
draw(Label("$\textbf{a}$", 0.2,W), l);

label("$E_2$", (1.5,-3,3.5));

import solids;
real r = 1;
currentlight=light(white, viewport=false,(5,-5,10));
draw(surface(sphere((0,0,0),r,100)),lightgray+opacity(0.4));
path3 c = Circle((0,0,0), 1, (1,-1,-1/2));
draw(c,black+linewidth(0.25mm));

import Projective2D;
metric=Metric(Elliptic);
MV P = Point(1,-3/5,4/5); //Point(1,-1/2,4/5); //Point(1,-1/2,1); 
MV a = join(Point(1,1,0),Point(1,0,2)); 

triple toxy(MV P) {return (1,topair(P).x,topair(P).y); }
dot(toxy(P)); label("$\textbf{P}$",toxy(P),(1/2,-1-1/2,1/2));
dot(toxy(dot(P,a)/a)); label("$(\textbf{P}\cdot\textbf{a})\textbf{a}^{-1}$",toxy(dot(P,a)/a),(2,-2,1.5));
dot(toxy(wedge(P,a)/a)); label("$(\textbf{P}\wedge\textbf{a})\textbf{a}^{-1}$",toxy(wedge(P,a)/a),(2,-2,1.5));

MV Q = dot(P,a)/a;
MV R = wedge(P,a)/a;
pen dot_pen = currentpen+2;
dot(normalise(toxy(P)), dot_pen); dot(-normalise(toxy(P)), dot_pen); draw(-normalise(toxy(P))--toxy(P));
dot(normalise(toxy(Q)), dot_pen); dot(-normalise(toxy(Q)), dot_pen); draw(-normalise(toxy(Q))--toxy(Q));
dot(normalise(toxy(R)), dot_pen); dot(-normalise(toxy(R)), dot_pen); draw(-normalise(toxy(R))--toxy(R));

path3 R = plane_path3(totriple(toline(dot(a,P))),(0,0,0),scale=10);
path3 lp = plane_intersection(R,p);
draw(Label("$\textbf{a}\cdot\textbf{P}$", 0.2,W), lp);

path3 c = Circle((0,0,0), 1, totriple(toline(dot(a,P))));
draw(c,black+linewidth(0.25mm));

//draw((1,-1/2,0)--(1,-1/2,1), Dotted);
//draw((1,0,1)--(1,-1/2,1), Dotted);

\end{asy}
\end{subfloatenv}\hspace{-12pt}%
\begin{subfloatenv}{\R{3}}
\begin{asy}
import Drawing3D;
DrawingR3 drawing = DrawingR3(4.0,0.1, camera=(10,2,7));	
drawing.target_axes();

path3 p = plane_path3((1,0,0),(1,0,0), v=(0,1,0));
drawing.Drawing3D.plane(p);
dot((1,0,0));
draw((1,-3.5,0)--(1,3.5,0));
draw((1,0,-3.5)--(1,0,3.5));

path3 Q = plane_path3((1,-1,-1/2), (1,1,0), theta=-13, scale=7.8);
path3 l = plane_intersection(Q,p);
draw(Label("$\textbf{a}$", 0.2,W), l);

import Projective2D;
metric=Metric(Elliptic);
MV a = join(Point(1,1,0),Point(1,0,2));
//MV b = join(Point(1,0,2),Point(1,1/6,0));
MV b = join(Point(1,0,0),Point(1,0,1));

path3 R = plane_path3(totriple(b), (0,0,0), scale=10);
path3 lb = plane_intersection(R,p);
draw(Label("$\textbf{b}$", 0.15,(1,-1,0)), lb);

label("$E_2$", (1.5,-3,3.5));

import solids;
real r = 1;
currentlight=light(white, viewport=false,(5,-5,10));
draw(surface(sphere((0,0,0),r,100)),lightgray+opacity(0.4));
path3 c = Circle((0,0,0), 1, (1,-1,-1/2));
draw(c,black+linewidth(0.25mm));

path3 C = plane_path3(totriple(-b*a/b), (0,0,0), scale=10);
path3 lc = plane_intersection(C,p);
draw(Label("$-\textbf{b}\textbf{a}\textbf{b}^{-1}$", 0.25,W), lc);

path3 cb = Circle((0,0,0), 1, totriple(b));
draw(cb,black+linewidth(0.25mm));

path3 cab = Circle((0,0,0), 1, totriple(-b*a/b));
draw(cab,black+linewidth(0.25mm));

\end{asy}
\end{subfloatenv}
\caption{Some basic properties of \El{2}}
\label{basic properties in El2}
\end{figure}

The line passing through \(\tb{P}\) and perpendicular to \(\tb{a}\) is given by 
the inner product \(\tb{a}\cdot\tb{P}\) as usual (see Figure~\ref{basic properties in El2}(a), where
\(\tb{P} = \e_{12}-\tfrac{3}{5}\e_{20}+\tfrac{4}{5}\e_{01}\) and \(\tb{a} = -2\e_0+2\e_1+\e_2\).
Note that \(\tb{a}\cdot\tb{P}=0\) if the location of \(\tb{P}\) 
coincides with that of the polar point \(\tb{a}\I\) of \(\tb{a}\). 
Projection, rejection, and reflection are defined in the usual way and since any non-zero blade is invertible in \El{2}
it is possible to project on the line \(\e_0\) and, say, the point \(\e_{20}\), for example.
For an illustration of the projection and rejection of a point by a line, see Figure~\ref{basic properties in El2}(a).
Note that the rejection \((\tb{P}\wedge\tb{a})\tb{a}^{-1}\) is at the distance of \(\tfrac{\pi}{2}\) from \(\tb{a}\),
so it is located at the same position in \El{2} as the polar point of \(\tb{a}\).
The top-down reflection of a line \(\tb{a}\) in a line \tb{b} is given by \(-\tb{b}\tb{a}\tb{b}^{-1}\)
and is illustrated in Figure~\ref{basic properties in El2}(b), where \(\tb{a} = -2\e_0+2\e_1+\e_2\) and \(\tb{b}=\e_1\).
In general, in \El{2} the projection of \(B_l\) on \(A_k\) is given by \((B_l\cdot A_k)A_k^{-1}\)
and the rejection of \(B_l\) by \(A_k\) is given by \((B_l\wedge A_k)A_k^{-1}\).
The top-down reflection of \(B_l\) in \(A_k\) is given by \((-1)^{kl}A_kB_lA_k^{-1}\)
and the bottom-up reflection by \(A_kB_lA_k^{-1}\).

Note that \(\tb{a}\cdot\tb{P}\) passes through the polar point \(\tb{a}\I\) of the line \(\tb{a}\).
In fact, any line in \El{2} passing through \(\tb{a}\I\) is perpendicular to \(\tb{a}\).
Moreover, the polar point of any line perpendicular to \(\tb{a}\) lies on the line \(\tb{a}\).
This can be used to construct a triangle with maximal area of \(\tfrac{\pi}{2}\) as follows.
For any line \(\tb{a}\), select a line \(\tb{b}\) passing through \(\tb{a}\I\), 
then define the third line \(\tb{c}=(\tb{a}\I)\vee(\tb{b}\I)\).
The angles of the resulting triangle are all equal to \(\tfrac{\pi}{2}\) and, therefore, its area equals \(\tfrac{\pi}{2}\).
Observe also that in elliptic space \El{2} any two lines intersect at a single point, which is at a finite distance from the origin.
Hence, there are no parallel lines in \El{2}.

\begin{figure}[t!]
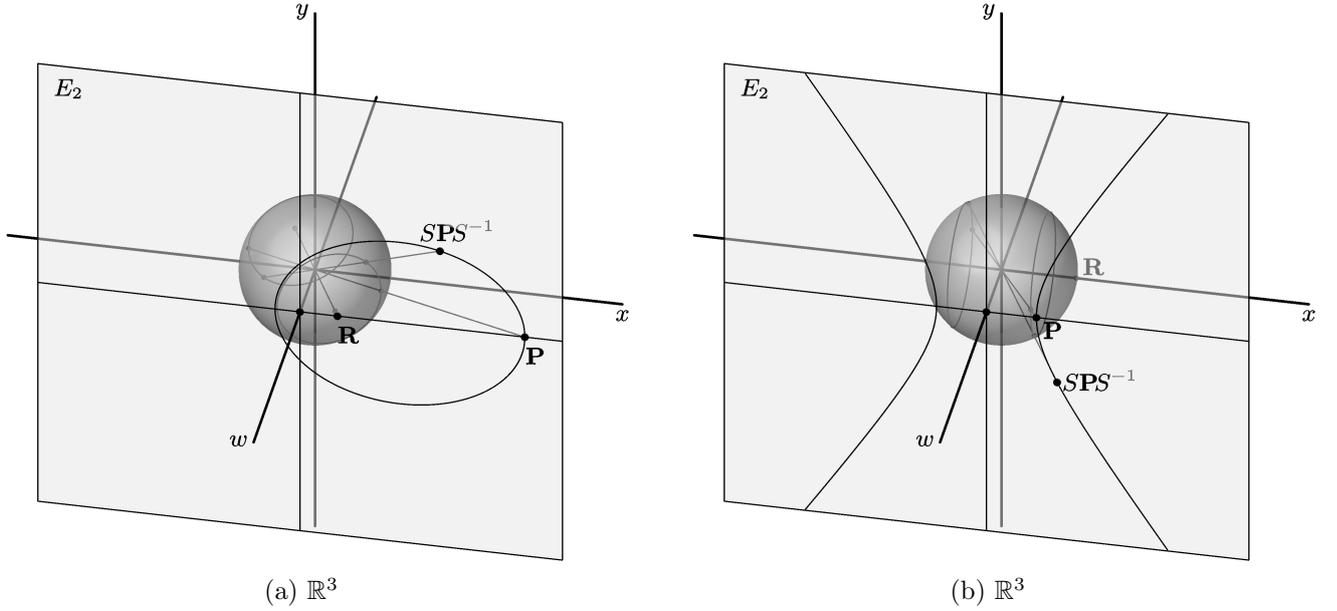

\begin{subfloatenv}{\R{3}}
\begin{asy}
import Drawing3D;
DrawingR3 drawing = DrawingR3(4.0,0.1, camera=(10,2,7));	
drawing.target_axes();

path3 p = plane_path3((1,0,0),(1,0,0), v=(0,1,0));
drawing.Drawing3D.plane(p);
dot((1,0,0));
draw((1,-3.5,0)--(1,3.5,0));
draw((1,0,-3.5)--(1,0,3.5));

label("$E_2$", (1.5,-3,3.5));

import solids;
real r = 1;
currentlight=light(white, viewport=false,(5,-5,10));
draw(surface(sphere((0,0,0),r,100)),lightgray+opacity(0.4));

import Projective2D;
metric=Metric(Elliptic);
MV R = Point(1,1/2,0); R/=norm(R);
MV P = Point(1,3,0); P/=norm(P);

real alpha = pi/4;
MV SS = cos(1/2*alpha) - R*sin(1/2*alpha); //exp(-1/2*alpha*R);
MV Q = SS*P/SS;

triple toxy(MV P) {return (1,topair(P).x,topair(P).y); };
pen dot_pen = currentpen+2;
dot(normalise(toxy(R)), dot_pen); dot(-normalise(toxy(R)), dot_pen);dot(toxy(R)); label("$\textbf{R}$",toxy(R),(1,1,-1));
dot(normalise(toxy(P)), dot_pen); dot(-normalise(toxy(P)), dot_pen); dot(toxy(P)); label("$\textbf{P}$",toxy(P),(1,1,-1));
dot(normalise(toxy(Q)), dot_pen); dot(-normalise(toxy(Q)), dot_pen); dot(toxy(Q)); label("$S\textbf{P}\!S^{-1}$",toxy(Q),(1,1,3));

draw(-normalise(toxy(R))--toxy(R));
draw(-normalise(toxy(P))--toxy(P));
draw(-normalise(toxy(Q))--toxy(Q));

triple to3(MV point) {return (point.w, point.x, point.y);};

triple f(real t) { MV S = cos(1/2*t)-R*sin(1/2*t); return to3(S*P/S); };
real x(real t) { return f(t).x; };
real y(real t) { return f(t).y; };
real z(real t) { return f(t).z; };
path3 c1 = graph(x, y, z, 0, 2*pi, operator ..);
draw(c1);

triple f(real t) { MV S = cos(1/2*t)-R*sin(1/2*t); return -to3(S*P/S); };
real x(real t) { return f(t).x; };
real y(real t) { return f(t).y; };
real z(real t) { return f(t).z; };
path3 c1 = graph(x, y, z, 0, 2*pi, operator ..);
draw(c1);

triple f(real t) { MV S = cos(1/2*t)-R*sin(1/2*t); return toxy(S*P/S); };
real x(real t) { return f(t).x; };
real y(real t) { return f(t).y; };
real z(real t) { return f(t).z; };
path3 c2 = graph(x, y, z, 0, 2*pi, operator ..);
draw(c2);

\end{asy}
\end{subfloatenv}\hspace{-12pt}%
\begin{subfloatenv}{\R{3}}
\begin{asy}
import Drawing3D;
DrawingR3 drawing = DrawingR3(4.0,0.1, camera=(10,2,7));	
drawing.target_axes();

path3 p = plane_path3((1,0,0),(1,0,0), v=(0,1,0));
drawing.Drawing3D.plane(p);
dot((1,0,0));
draw((1,-3.5,0)--(1,3.5,0));
draw((1,0,-3.5)--(1,0,3.5));

label("$E_2$", (1.5,-3,3.5));

import solids;
real r = 1;
currentlight=light(white, viewport=false,(5,-5,10));
draw(surface(sphere((0,0,0),r,100)),lightgray+opacity(0.4));

import Projective2D;
metric=Metric(Elliptic);
MV R = Point(0,1,0); R/=norm(R);
MV P = Point(1,2/3,0); P/=norm(P);

real alpha = pi/4;
MV SS = cos(1/2*alpha) - R*sin(1/2*alpha); //exp(-1/2*alpha*R);
MV Q = SS*P/SS;

triple toxy(MV P) {return (1,topair(P).x,topair(P).y); };
pen dot_pen = currentpen+2;
dot(normalise(toxy(P)), dot_pen); dot(-normalise(toxy(P)), dot_pen); dot(toxy(P)); label("$\textbf{P}$",toxy(P),(1,1.25,-0.2));
dot(normalise(toxy(Q)), dot_pen); dot(-normalise(toxy(Q)), dot_pen); dot(toxy(Q)); label("$S\textbf{P}\!S^{-1}$",toxy(Q),(1.5,2.25,1.5));

triple to3(MV point) {return (point.w, point.x, point.y);};
dot(to3(R), dot_pen); dot(-to3(R), dot_pen); label("$\textbf{R}$",to3(R),(0,1,1));

draw(-normalise(toxy(P))--toxy(P));
draw(-normalise(toxy(Q))--toxy(Q));

triple to3(MV point) {return (point.w, point.x, point.y);};

triple f(real t) { MV S = cos(1/2*t)-R*sin(1/2*t); return to3(S*P/S); };
real x(real t) { return f(t).x; };
real y(real t) { return f(t).y; };
real z(real t) { return f(t).z; };
path3 c1 = graph(x, y, z, 0, 2*pi, operator ..);
draw(c1);

triple f(real t) { MV S = cos(1/2*t)-R*sin(1/2*t); return -to3(S*P/S); };
real x(real t) { return f(t).x; };
real y(real t) { return f(t).y; };
real z(real t) { return f(t).z; };
path3 c1 = graph(x, y, z, 0, 2*pi, operator ..);
draw(c1);

triple f(real t) { MV S = cos(1/2*t)-R*sin(1/2*t); return toxy(S*P/S); };
real x(real t) { return f(t).x; };
real y(real t) { return f(t).y; };
real z(real t) { return f(t).z; };
path3 c2 = graph(x, y, z, -1.291, 1.291, operator ..);
draw(c2);

path3 c2 = graph(x, y, z, pi-1.291, pi, operator ..);
draw(c2);

path3 c2 = graph(x, y, z, -pi, -pi+1.291, operator ..);
draw(c2);

\end{asy}
\end{subfloatenv}
\caption{Rotation in \El{2}}
\label{rotation in El2}
\end{figure}

The spin group of \El{2} is defined in the standard way and consists of spinors of the form \(e^A\) where \(A\) is a bivector.
It is isomorphic to the group of normalised quaternions.
The action \(SA_kS^{-1}\) of a spinor \(S=e^{-\tfrac{1}{2}\alpha\tb{R}}\) 
on the geometric object dually represented  by the blade \(A_k\)
is consistent with the rotation of \(A_k\) around the point \(\tb{R}\) by the angle \(\alpha\),
where \(\alpha\) refers to the angle according to the elliptic metric, rather than Euclidean.
No translation is possible in \El{2} since there are no points at infinity.

The rotation of the point \(\tb{P}=\e_{12}+3\e_{20}\) around \(\tb{R}=\e_{12}+\tfrac{1}{2}\e_{20}\)
by the angle \(\alpha=\tfrac{\pi}{4}\) is shown in Figure~\ref{rotation in El2}(a).
The solid curve in \El{2} (see Figure~\ref{rotation in El2}(a)) shows the trajectory of \(\tb{P}\)
under the action of the spinor \(S=e^{-\tfrac{1}{2}t\tb{R}}\), where the parameter \(t\)
varies from \(0\) to \(2\pi\) (the projection of this curve onto the unit sphere is shown with the solid curve as well).
The curve consists of points at an equal distance from \(\tb{R}\) and, therefore, forms a circle in \El{2},
which passes through \(\tb{P}\) and whose centre is at \(\tb{R}\).
The appearance of the circle on the plane \El{2} is elliptic.

Figure~\ref{rotation in El2}(b) shows another example of the rotation, where 
\(\tb{P}=\e_{12}+\tfrac{2}{3}\e_{20}\), \(\tb{R}=\e_{20}\), and \(\alpha=\tfrac{\pi}{4}\). 
Note that \(\tb{R}\) is shown on the unit sphere only; it corresponds to a stack of lines in the plane \El{2}.
The solid curve shows the trajectory of \(\tb{P}\) 
under the action of the spinor \(S=e^{-\tfrac{1}{2}t\tb{R}}\), where the parameter \(t\) varies from \(0\) to \(2\pi\).
It is also a circle but its appearance in the plane \El{2} is hyperbolic.

In general, a circle in \El{2} could appear as an ellipse, two branches of a hyperbola, a parabola, 
or a straight line.
The latter represents a circle  whose radius equals \(\tfrac{\pi}{2}\).
The circle is parabolic if there is only one point on it that belongs to the line \(\e_0\)
and hyperbolic if there are two points on the circle that belong to \(\e_0\).
Otherwise, it is elliptic.
The circle of radius \(\tfrac{\pi}{2}\) whose centre is at the origin cannot be visualised on the plane \El{2};
this circle corresponds to the line \(\e_0\).

\section{Elliptic space \El{3}}
For a plane \(\tb{a}=d\e_0+a\e_1+b\e_2+c\e_3\), 
a line \(\mb{\Lambda}=p_{10}\e_{10}+p_{20}\e_{20}+p_{30}\e_{30}+p_{23}\e_{23}+p_{31}\e_{31}+p_{12}\e_{12}\),
where \(p_{10}p_{23}+p_{20}p_{31}+ p_{30}p_{12}=0\), and 
a point \(\tb{P}=w\e_{123}+x\e_{320}+y\e_{130}+z\e_{210}\), the norm is given by 
\(\norm{\tb{a}}=\sqrt{d^2+a^2+b^2+c^2}\),
\(\norm{\mb{\Lambda}}=\sqrt{p_{10}^2+p_{20}^2+p_{30}^2+p_{23}^2+p_{31}^2+p_{12}^2}\),
\(\norm{\tb{P}}=\sqrt{w^2+x^2+y^2+z^2}\),
and \(\tb{a}^2=\norm{\tb{a}}^2\), \(\mb{\Lambda}^2=-\norm{\mb{\Lambda}}^2\), \(\tb{P}^2=-\norm{\tb{P}}^2\).
Note that the above formula for \(\norm{\mb{\Lambda}}\) is only valid for lines; it is not applicable to non-simple bivectors.

The hyperplane \El{3} is embedded in the model space \R{4} at \(w=1\).
By analogy with the 1- and 2-dimensional cased, the 3-dimensional elliptic space 
can also be thought of as the unit 3-sphere centred on the origin of \R{4}, which is defined by
\begin{equation}
w^2+x^2+y^2+z^2=1.
\end{equation}
The antipodal points in the unit 3-sphere are identified.
The immediate 3-dimensional vicinity of the origin of \El{3} is  near the unit 3-sphere
and they resemble each other closely. 
As one moves away from the origin of \El{3}, the unit 3-sphere deviates
from the hyperplane \El{3}, curving into the fourth dimension along \(w\)
and towards the hyperplane \(w=0\).
It is not possible to visualise the unit 3-sphere.
However, one can infer the geometry in the 3-sphere by examining the geometry in \El{3},
which can be readily visualised.

Planes and lines in \El{3} correspond to great spheres and great circles in the unit 3-sphere.
For instance, the plane \(\e_0\) corresponds to the great sphere
 where the unit 3-sphere intersects the hyperplane \(w=0\).
Any two distinct great spheres in the unit 3-sphere intersect along a great circle,
which corresponds to the line at the intersect of two planes in \El{3} corresponding to the great spheres.
A point in \El{3} corresponds to a pair of antipodal points in the unit 3-sphere.

\begin{figure}[t]
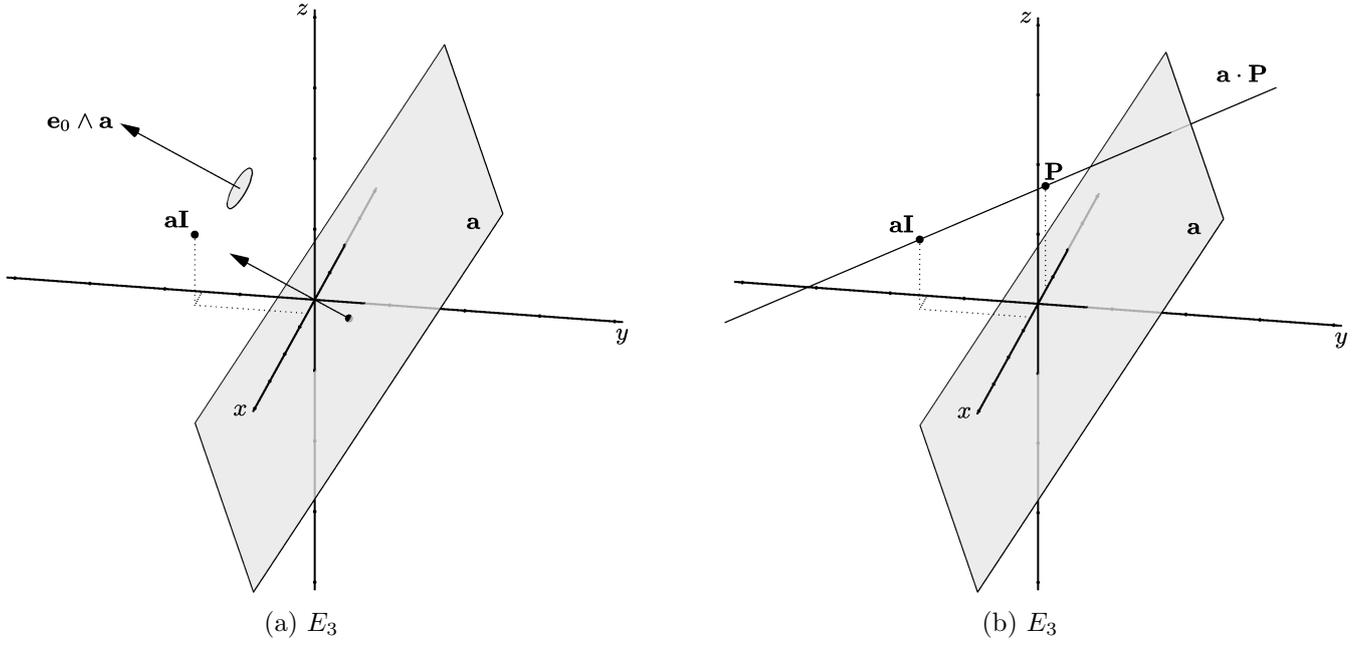
\hspace{-1cm}
\begin{subfloatenv}{\El{3}}
\begin{asy}
import Figure3D;
Figure f = Figure();
metric=Metric(Elliptic);

var a = Plane(1,1/2,-3/2,1);
f.plane(a, "$\textbf{a}$", align=(1,-3,-1));

var L = wedge(e_0,a);
f.line_at_infinity(L, centre=(0,-1,1.5), label="$\textbf{e}_0\wedge\textbf{a}$", align=(0,-1,0));

var N = a*I;
f.point(N, "$\textbf{a}\textbf{I}$", draw_orientation=false, align=(0,-1,1));

\end{asy}
\end{subfloatenv}\hfill%
\begin{subfloatenv}{\El{3}}
\begin{asy}
import Figure3D;
Figure f = Figure();
metric=Metric(Elliptic);

var a = Plane(1,1/2,-3/2,1);
f.plane(a, "$\textbf{a}$", align=(1,-3,-1),draw_orientation=false);

var P = Point(1,-1/2,0,1.5); //Point(1,-1,1,1.5);
f.point(P, "$\textbf{P}$", draw_orientation=false, align=(0,0.5,1));

var N = a*I;
f.point(N, "$\textbf{a}\textbf{I}$", draw_orientation=false, align=(0,-1,1));

var L = dot(a,P);
f.line(L, "$\textbf{a}\cdot\textbf{P}$", position=0.02, draw_orientation=false, align=(0,-1,1.5));

\end{asy}
\end{subfloatenv}
\caption{Basic properties of points and planes  in \El{3}}
\label{basic  E3}
\end{figure}

The polar point \(\tb{a}\I\) of a plane \(\tb{a}\) is shown in Figure~\ref{basic  E3}(a),
where \(\tb{a}=\e_0+\tfrac{1}{2}\e_1-\tfrac{3}{2}\e_2+\e_3\). 
The line \(\e_0\wedge\tb{a}\) is at the intersection of the planes \(\tb{a}\) and \(\e_0\).
It is at a finite distance from the origin but it cannot be shown directly since it lies in \(\e_0\)
(it is depicted in the same way as a line at infinity in \E{3}).
A line passing through \(\tb{P}\) and perpendicular to \(\tb{a}\) is given by \(\tb{a}\cdot\tb{P}\),
unless \(\tb{P}\) is located at \(\tb{a}\I\) in which case \(\tb{a}\cdot\tb{P}=0\).
An example is shown in Figure~\ref{basic  E3}(b) where \(\tb{P}=\e_{123}-\tfrac{1}{2}\e_{320}+\tfrac{3}{2}\e_{210}\).
Note that \(\tb{a}\cdot\tb{P}\) passes through the polar point \(\tb{a}\I\) of \(\tb{a}\).
In fact, any line perpendicular to \(\tb{a}\) passes through \(\tb{a}\I\) and
and conversely any line passing through \(\tb{a}\I\) is perpendicular to the plane  \(\tb{a}\).
The same observation applies to planes, i.e.\ a plane is perpendicular to \(\tb{a}\) if and only if
it passes through \(\tb{a}\I\).
This also concerns the plane passing through \(\tb{a}\I\) which would be parallel to \(\tb{a}\)
in Euclidean space (in elliptic space, there are no parallel planes since any two planes intersect
at a finite distance from the origin). 
Note also that any line \(\mb{\Lambda}\) that lies in the plane \(\tb{a}\) is in the same relationship with respect to 
\(\tb{a}\I\) as an equatorial line and its polar point in \El{2}.

The distance \(r\in[0,\tfrac{\pi}{2}]\) between normalised points \(\tb{P}\) and \(\tb{Q}\) is defined by
\begin{equation}
\sin r= \norm{\tb{P}\vee\tb{Q}},
\end{equation}
where \(\tb{P}\vee\tb{Q}\) is a line passing through \(\tb{P}\) and \(\tb{Q}\).
Since \(\norm{\tb{P}\vee\tb{Q}}^2+|\tb{P}\cdot\tb{Q}|^2=1\) in \El{3}, 
one also gets \(\cos r= |\tb{P}\cdot\tb{Q}|\).
If \(\tb{P}=\e_{123}+x\e_{320}+y\e_{130}+z\e_{210}\) and \(\tb{Q}=\e_{123}\), i.e.\ \(\tb{Q}\)
is at the origin of \El{3}, then
\(\norm{\tb{P}}=\sqrt{1+x^2+y^2+z^2}\) and the distance from the origin to \(\tb{P}\)
satisfies
\begin{equation}
\sin r = \frac{\sqrt{x^2+y^2+z^2}}{\sqrt{1+x^2+y^2+z^2}}, \quad \cos r = \frac{1}{\sqrt{1+x^2+y^2+z^2}}.
\end{equation}
As in 1- and 2-dimensional cases, the maximum distance in \El{3} equals \(\tfrac{\pi}{2}\).

\begin{figure}[t]
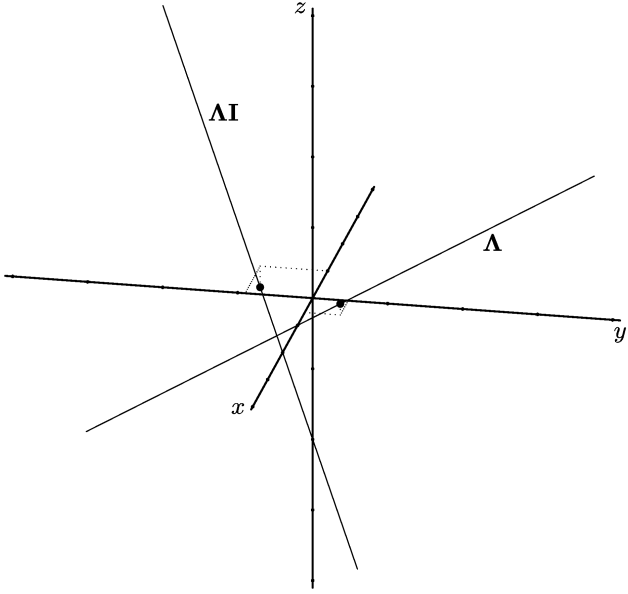
\hspace{-1cm}
\begin{asy}
import Figure3D;
Figure f = Figure();
metric=Metric(Elliptic);

var P = Point(1,1,0,0);
var Q = Point(1,0,1,1/3);
var L = join(P,Q);
write(L);

f.line(L, "$\boldsymbol{\Lambda}$", align=(0,0,-1),draw_orientation=false);

MV C = toline(L).centre();
f.point(C," ",draw_orientation=false,draw_helper_lines=true);

f.line(L*I, "$\boldsymbol{\Lambda}\textbf{I}$", position=0.8,align=(0,1,0.5),draw_orientation=false);

MV CI = toline(L*I).centre();
f.point(CI," ",draw_orientation=false,draw_helper_lines=true);

\end{asy}
\caption{Lines in \El{3}}
\label{lines in El3}
\end{figure}

A plane \(\tb{a}\) and its polar point \(\tb{a}\I\) are related by projective duality
and \(\J(\tb{a})=\Id(\tb{a}\I)\) where \(\J\) and \(\Id\) are duality and identity transformations
(note that \(\I^{-1}=\I\) and \(\I^2=1\) in \El{3}).
Furthermore, the distance between \(\tb{a}\I\) and any point on \(\tb{a}\) equals \(\tfrac{\pi}{2}\).
The polar point of the plane \(\tb{P}\I\) coincides with \(\tb{P}\) and \(\J(\tb{P})=\Id(\tb{P}\I)\).
The distance from \(\tb{P}\) to any point on the plane \(\tb{P}\I\) equals \(\tfrac{\pi}{2}\).
For a line \(\mb{\Lambda}\), the pair of lines \(\mb{\Lambda}\) and \(\mb{\Lambda}\I\) are
also related by projective duality and \(\J(\mb{\Lambda})=\Id(\mb{\Lambda}\I)\).
An example is shown in Figure~\ref{lines in El3}(a) 
where \(\mb{\Lambda}=-\tfrac{1}{3}\e_{20}+\e_{30}+\e_{23}-\e_{31}-\tfrac{1}{3}\e_{12}\).
The distance between any point on  \(\mb{\Lambda}\) and any point on \(\mb{\Lambda}\I\)
equals \(\tfrac{\pi}{2}\).

The angle \(\alpha\) between normalised planes \(\tb{a}\) and \(\tb{b}\),
or more precisely the angle between the top-down orientation vectors of the planes, is defined by
\begin{equation}
\cos\alpha=\tb{a}\cdot\tb{b}.
\end{equation}
For any normalised \(\tb{a}\) and \(\tb{b}\), their geometric product yields \(\tb{a}\tb{b}=\cos\alpha+\mb{\Lambda}\sin\alpha\),
where \(\mb{\Lambda}\) is a normalised line at the intersection of the planes \(\tb{a}\) and \(\tb{b}\).

A plane in \El{3} can be considered as a copy of \El{2}, so the total area of any plane 
in \El{3} is finite and is equal to \(2\pi\).
The total volume of \El{3} equals \(\pi^2\).
This is equivalent to half the volume of the unit 3-sphere in \R{4} with the Euclidean metric. 
A factor of \(\tfrac{1}{2}\) is needed since antipodal points in the unit 3-sphere serving as a model of \El{3} are identified.
The volume of the unit 3-sphere in \R{4} should not be confused
with the hypervolume of the 4-dimensional region bounded by the unit 3-sphere.
There is no simple expression for the volume of an arbitrary 3-simplex, i.e.\ a tetrahedron, in \El{3}.
However, the volume of a 3-simplex whose sides are perpendicular to one another equals \(\tfrac{\pi^2}{8}\),
which is the largest possible volume a 3-simplex can have in \El{3}.

The distance \(r\) from a normalised point \(\tb{P}\) to a normalised plane \(\tb{a}\) satisfies \(\sin r=|\tb{a}\vee\tb{P}|\)
and \(\cos r = \norm{\tb{a}\cdot\tb{P}}\).
It equals the distance along the line \(\tb{a}\cdot\tb{P}\) from \(\tb{P}\) to the point where \(\tb{a}\cdot\tb{P}\) intersects \(\tb{a}\).
The distance from a normalised point \(\tb{P}\) 
to a normalised line \(\mb{\Lambda}\) satisfies \(\sin r=\norm{\mb{\Lambda}\vee\tb{P}}\)
and \(\cos r = \norm{\mb{\Lambda}\cdot\tb{P}}\), 
where \(\mb{\Lambda}\vee\tb{P}\) is a plane passing through \(\mb{\Lambda}\) and \(\tb{P}\),
and \(\mb{\Lambda}\cdot\tb{P}\) is a plane perpendicular to \(\mb{\Lambda}\) and passing through \(\tb{P}\).
Recall that \(\mb{\Lambda}\vee\tb{P}=0\) if \(\tb{P}\) lies on the line \(\mb{\Lambda}\)
and note that \(\mb{\Lambda}\cdot\tb{P}=0\) if \(\tb{P}\) lies on \(\mb{\Lambda}\I\).
The angle \(\alpha\) between a line \(\mb{\Lambda}\) and a plane \(\tb{a}\) 
satisfies \(\cos\alpha=\norm{\tb{a}\cdot\mb{\Lambda}}\) and \(\sin\alpha=\norm{\tb{a}\wedge\mb{\Lambda}}\),
where \(\tb{a}\cdot\mb{\Lambda}\) is a plane passing through \(\mb{\Lambda}\) and perpendicular to \(\tb{a}\),
and \(\tb{a}\wedge\mb{\Lambda}\) is a point where \(\mb{\Lambda}\) intersects \(\tb{a}\)
(the angle \(\alpha\) does not take into account the orientation of \(\mb{\Lambda}\) and \(\tb{a}\)).
Recall that \(\tb{a}\wedge\mb{\Lambda}=0\) if the line \(\mb{\Lambda}\) lies in the plane \(\tb{a}\).
Note also that \(\tb{a}\cdot\mb{\Lambda}=0\) if the line \(\mb{\Lambda}\) passes through \(\tb{a}\I\),
in which case \(\mb{\Lambda}\I\) lies in the plane \(\tb{a}\).
A plane is perpendicular to \(\mb{\Lambda}\) if and only if it passes through \(\mb{\Lambda}\I\).
So, a set of planes perpendicular to \(\mb{\Lambda}\) forms a sheaf of planes 
passing through  \(\mb{\Lambda}\I\).
A line is perpendicular to a plane \(\tb{a}\) if and only if it passes through \(\tb{a}\I\).
The angle \(\alpha\) between two normalised lines \(\mb{\Lambda}\) and \(\mb{\Theta}\) is given by
\(\mb{\Lambda}\cdot\mb{\Theta}=-\cos\alpha\) if the lines intersect.


The pair of lines \(\mb{\Lambda}\) and \(\mb{\Lambda}\I\) have an interesting property that
the distance between \(\mb{\Lambda}\) and \(\mb{\Lambda}\I\) is constant everywhere along these lines.
Such lines in \El{3} are called Clifford-parallel.
There are two families, positive and negative, of Clifford-parallel lines associated with the line \(\mb{\Lambda}\).
The lines in each family are Clifford-parallel to one another and \(\mb{\Lambda}\)
is Clifford-parallel to lines in both families.
Each family foliates the whole elliptic space \El{3}. 
If one ignores weight and orientation,
there are exactly two lines Clifford-parallel to \(\mb{\Lambda}\)
passing through every point in \El{3}, except for the points on the lines \(\mb{\Lambda}\) and \(\mb{\Lambda}\I\).
One of these parallels is positive and the other is negative. 

\begin{figure}[t!]
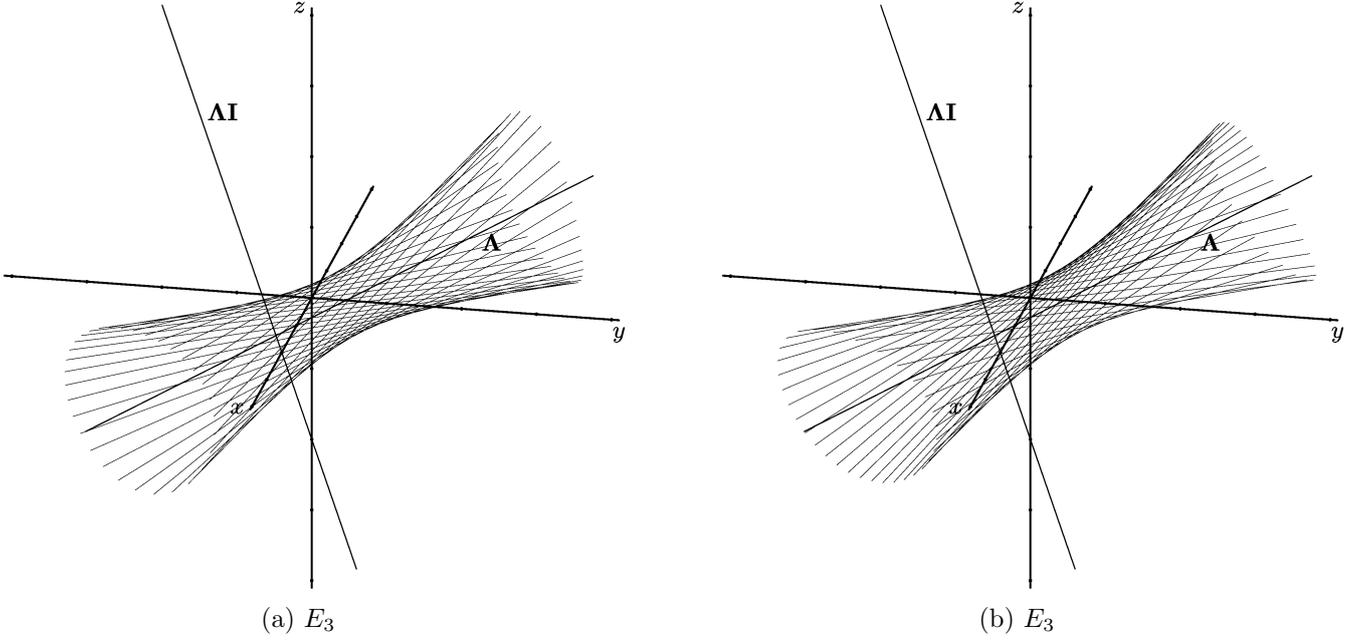
\hspace{-1cm}
\begin{subfloatenv}{\El{3}}
\begin{asy}
import Figure3D;
Figure f = Figure();
metric=Metric(Elliptic);

var P = Point(1,1,0,0);
var Q = Point(1,0,1,1/3);
var L = join(P,Q);
write(L);

f.line(L, "$\boldsymbol{\Lambda}$", align=(0,0,-1),draw_orientation=false);

f.line(L*I, "$\boldsymbol{\Lambda}\textbf{I}$", position=0.8,align=(0,1,0.5),draw_orientation=false);

L/=norm(L);
Clifford_parallel c = Clifford_parallel(L);
real theta=pi/2-pi/10;
int N=32;
for(int i: sequence(N)) { real phi = 2*pi*i/N; f.line(c.positive(phi,theta), "", draw_orientation=false,size=3.5,pen=currentpen+0.25);};

\end{asy}
\end{subfloatenv}\hfill%
\begin{subfloatenv}{\El{3}}
\begin{asy}
import Figure3D;
Figure f = Figure();
metric=Metric(Elliptic);

var P = Point(1,1,0,0);
var Q = Point(1,0,1,1/3);
var L = join(P,Q);
write(L);

f.line(L, "$\boldsymbol{\Lambda}$", align=(0,0,-1),draw_orientation=false);

f.line(L*I, "$\boldsymbol{\Lambda}\textbf{I}$", position=0.8,align=(0,1,0.5),draw_orientation=false);

L/=norm(L);
Clifford_parallel c = Clifford_parallel(L);
real theta=pi/2-pi/10;
int N=32;
for(int i: sequence(N)) { real phi = 2*pi*i/N;  f.line(c.negative(phi,theta), "", draw_orientation=false,size=3.5,pen=currentpen+0.25);};

\end{asy}
\end{subfloatenv}
\caption{Clifford parallels of \(\mb{\Lambda}\)}
\label{Clifford parallels in El3}
\end{figure}

I assume the line \(\mb{\Lambda}\) is normalised and is given by
\(\mb{\Lambda}=p_{10}\e_{10}+p_{20}\e_{20}+p_{30}\e_{30}+p_{23}\e_{23}+p_{31}\e_{31}+p_{12}\e_{12}\)
with \(p_{10}p_{23}+p_{20}p_{31}+ p_{30}p_{12}=0\).
The positive family of lines Clifford-parallel to \(\mb{\Lambda}\)  is given by
\begin{equation}
\mb{\Lambda}^\p=\mb{\Lambda}^\p(\mb{\Omega})=\mb{\Lambda}+(\mb{\Lambda}_m\cdot\mb{\Omega})(\I-1)\mb{\Omega},
\label{positive Clifford parallel}
\end{equation}
and the negative family by
\begin{equation}
\mb{\Lambda}^\n=\mb{\Lambda}^\n(\mb{\Omega})=\mb{\Lambda}+(\mb{\Lambda}_p\cdot\mb{\Omega})(\I+1)\mb{\Omega},
\label{negative Clifford parallel}
\end{equation}
where 
\begin{equation}
\mb{\Lambda}_m=(p_{10}-p_{23})\e_{23}+(p_{20}-p_{31})\e_{31}+(p_{30}-p_{12})\e_{12},
\end{equation}%
\begin{equation}
\mb{\Lambda}_p=(p_{10}+p_{23})\e_{23}+(p_{20}+p_{31})\e_{31}+(p_{30}+p_{12})\e_{12},
\end{equation}
and \(\mb{\Omega}\) ranges over all normalised lines passing through the origin.
The lines \(\mb{\Lambda}_m\) and \(\mb{\Lambda}_p\) also pass through the origin and are normalised since \(\mb{\Lambda}\) 
is assumed to be normalised.
The parallels \(\mb{\Lambda}^\p\), \(\mb{\Lambda}^\n\) are all normalised and their orientation is consistent with that of \(\mb{\Lambda}\).
It is convenient to parametrise \(\mb{\Omega}\) with \(\phi\in[0,2\pi)\) and \(\theta\in[0,\pi]\) by
\begin{equation}
\mb{\Omega}=\mb{\Omega}_m(\phi,\theta)=
e^{-\tfrac{1}{2}\phi\mb{\Lambda}_m}e^{-\tfrac{1}{2}\theta\mb{\Lambda}_m^\perp}
\mb{\Lambda}_m e^{\tfrac{1}{2}\theta\mb{\Lambda}_m^\perp}e^{\tfrac{1}{2}\phi\mb{\Lambda}_m}
\label{Omega m}
\end{equation}
for the positive family and by 
\begin{equation}
\mb{\Omega}=\mb{\Omega}_p(\phi,\theta)=
e^{-\tfrac{1}{2}\phi\mb{\Lambda}_p}e^{-\tfrac{1}{2}\theta\mb{\Lambda}_p^\perp}
\mb{\Lambda}_p e^{\tfrac{1}{2}\theta\mb{\Lambda}_p^\perp}e^{\tfrac{1}{2}\phi\mb{\Lambda}_p}
\label{Omega p}
\end{equation}
for the negative family;
\(\mb{\Lambda}_m^\perp\) is perpendicular to \(\mb{\Lambda}_m\),
 \(\mb{\Lambda}_p^\perp\) is perpendicular to \(\mb{\Lambda}_p\),
and both \(\mb{\Lambda}_m^\perp\) and \(\mb{\Lambda}_p^\perp\) are normalised and pass through the origin.
Then, \(|\tfrac{\pi}{2}-\theta|\) gives the distance from \(\mb{\Lambda}\)
to the Clifford parallels, both positive and negative, parametrised  with \((\phi,\theta)\).
The exact meaning of \(\phi\) depends on the choice of \(\mb{\Lambda}_m^\perp\) and  \(\mb{\Lambda}_p^\perp\).
Note that \(\mb{\Lambda}_m\cdot\mb{\Omega}_m(\phi,\theta)=-\cos\theta\), i.e.\ 
\(\theta\) in (\ref{Omega m}) is the angle between \(\mb{\Lambda}_m\) and \(\mb{\Omega}_m(\phi,\theta)\),
and  \(\mb{\Lambda}_p\cdot\mb{\Omega}_p(\phi,\theta)=-\cos\theta\), i.e.\ 
\(\theta\) in (\ref{Omega p}) is the angle between \(\mb{\Lambda}_p\) and \(\mb{\Omega}_p(\phi,\theta)\).
 So, the Clifford parallels can be given by
\begin{equation}
\mb{\Lambda}^\p=\mb{\Lambda}^\p(\phi,\theta)=\mb{\Lambda}-\cos\theta(\I-1)\mb{\Omega}_m(\phi,\theta),
\label{positive Clifford parallel parametrised}
\end{equation}
\begin{equation}
\mb{\Lambda}^\n=\mb{\Lambda}^\n(\phi,\theta)=\mb{\Lambda}-\cos\theta(\I+1)\mb{\Omega}_p(\phi,\theta).
\label{negative Clifford parallel parametrised}
\end{equation}
If \(\mb{\Lambda}\) is not normalised, one can find the Clifford parallels of \(\mb{\Lambda}/\norm{\mb{\Lambda}}\)
by applying the above formulas and then scale the parallels by \(\norm{\mb{\Lambda}}\) to get the
Clifford parallels of \(\mb{\Lambda}\).
Some positive and negative Clifford parallels of \(\mb{\Lambda}=\tfrac{3}{\sqrt{29}}(-\tfrac{1}{3}\e_{20}+\e_{30}+\e_{23}-\e_{31}-\tfrac{1}{3}\e_{12})\) 
located at the distance \(r=\tfrac{\pi}{10}\) from \(\mb{\Lambda}\)
 are shown in Figure~\ref{Clifford parallels in El3} (positive in (a) and negative in (b)).

Note that \(\mb{\Lambda}\I\) is a positive Clifford parallel of \(\mb{\Lambda}\)
(substitute  \(\mb{\Omega}=\mb{\Lambda}_m\) into  (\ref{positive Clifford parallel}))
and \(-\mb{\Lambda}\I\) is a negative Clifford parallel of \(\mb{\Lambda}\)
(substitute  \(\mb{\Omega}=\mb{\Lambda}_p\) into  (\ref{negative Clifford parallel})),
while \(\mb{\Lambda}\) belongs to both positive and negative families.
Also,  even though \(\mb{\Lambda}_m\) is used in the definition of the positive Clifford parallels,
 \(-\mb{\Lambda}_m\) is a negative Clifford parallel and \(\mb{\Lambda}_p\) is a positive Clifford parallel.
They can be obtained by substituting 
\(\mb{\Omega}=(\mb{\Lambda}_p+\mb{\Lambda}_m)/\norm{\mb{\Lambda}_p+\mb{\Lambda}_m}\)
into (\ref{positive Clifford parallel}) and (\ref{negative Clifford parallel}).

Since  \((\I+1)(\I-1)=0\), the positive and negative Clifford parallels of \(\mb{\Lambda}\) satisfy
\begin{equation}
\textstyle
(\I+1)\mb{\Lambda}^\p=(\I+1)\mb{\Lambda}\quad\text{and}\quad
(\I-1)\mb{\Lambda}^\n=(\I-1)\mb{\Lambda},
\label{parallels property}
\end{equation}
respectively.
These formulas can be used to define Clifford parallels.
A non-simple bivector which can be written as \(\mb{\Xi}^+=(\I+1)\mb{\Lambda}\) where \(\mb{\Lambda}\) is a line
has the property that \(\I\mb{\Xi}^+=\mb{\Xi}^+\). I will call it a positive Clifford bivector.
It defines a set of positive Clifford parallels via \((\I+1)\mb{\Lambda}^\p=\mb{\Xi}^+\), which includes 
\(\mb{\Lambda}\) and \(\mb{\Lambda}\I\).
I will refer to these positive Clifford parallels as the Clifford parallels of the positive Clifford bivector \(\mb{\Xi}^+\).
If \(\mb{\Lambda}^\p\) is a Clifford parallel of \(\mb{\Xi}^+\), then \(\mb{\Lambda}^\p\I\) is also a Clifford parallel of \(\mb{\Xi}^+\).
Similarly, a non-simple bivector which can be written as \(\mb{\Xi}^-=(\I-1)\mb{\Lambda}\) where \(\mb{\Lambda}\) is a line
has the property that \(\I\mb{\Xi}^-=-\mb{\Xi}^-\) and will be called a negative Clifford bivector.
It defines a set of negative Clifford parallels via \((\I-1)\mb{\Lambda}^\n=\mb{\Xi}^-\), which includes 
\(\mb{\Lambda}\) and \(-\mb{\Lambda}\I\).
I will refer to these negative Clifford parallels as the Clifford parallels of the negative Clifford bivector \(\mb{\Xi}^-\).
If \(\mb{\Lambda}^\n\) is a Clifford parallel of \(\mb{\Xi}^-\), then \(-\mb{\Lambda}^\n\I\) is also a Clifford parallel of \(\mb{\Xi}^-\).
The Clifford parallel \(\mb{\Lambda}\) of a Clifford bivector \(\mb{\Xi}\) which passes through a point \(\tb{P}\) is given by 
\begin{equation}
\mb{\Lambda}=(\mb{\Xi}\vee\tb{P})\tb{P}^{-1}.
\end{equation}
It is either positive or negative depending on whether the Clifford bivector \(\mb{\Xi}\) is positive or negative.

Most bivectors in \El{3} are non-simple, but any non-simple bivector can be decomposed into a sum of two complementary
lines called the axes of the non-simple bivector.
Two lines \(\mb{\Lambda}_1\) and \(\mb{\Lambda}_2\) are called complementary if 
\(\mb{\Lambda}_1\cdot\mb{\Lambda}_2=0\) and \(\mb{\Lambda}_1\times\mb{\Lambda}_2=0\).
Complementary lines are not only perpendicular but they also commute.
For instance, \(\mb{\Theta}\) and \(\mb{\Theta}\I\), where \(\mb{\Theta}\) is a line, are complementary.
For a non-simple \(\mb{\Lambda}\), assume \(\mb{\Lambda}=\mb{\Lambda}_1+\mb{\Lambda}_2\), where
\(\mb{\Lambda}_1\) and \(\mb{\Lambda}_2\) are complementary.
Then
\begin{equation}
\mb{\Lambda}_1\mb{\Lambda}_2=\tfrac{1}{2}\mb{\Lambda}\wedge\mb{\Lambda},
\end{equation}
\begin{equation}
\mb{\Lambda}_1^2+\mb{\Lambda}_2^2=\mb{\Lambda}\cdot\mb{\Lambda}.
\end{equation}
Squaring the first equation and using the second yields a quadratic equation 
for \(\mb{\Lambda}_1^2\) and \(\mb{\Lambda}_2^2\):
\begin{equation}
\mb{\Lambda}_{1,2}^4-(\mb{\Lambda}\cdot\mb{\Lambda})\mb{\Lambda}_{1,2}^2
+\tfrac{1}{4}(\mb{\Lambda}\vee\mb{\Lambda})^2=0.
\end{equation}
If \((\mb{\Lambda}\cdot\mb{\Lambda})^2>(\mb{\Lambda}\vee\mb{\Lambda})^2\), 
there are two distinct solutions given by
\begin{equation}
\mb{\Lambda}_{1,2}^2=\tfrac{1}{2}\left(\mb{\Lambda}\cdot\mb{\Lambda}
\pm\sqrt{(\mb{\Lambda}\cdot\mb{\Lambda})^2-(\mb{\Lambda}\vee\mb{\Lambda})^2}\right).
\label{L1^2 and L2^2}
\end{equation}
For definiteness, I will assume \(\mb{\Lambda}_1\) uses the minus sign in (\ref{L1^2 and L2^2})
and \(\mb{\Lambda}_2\) uses the plus sign, which implies 
\(\norm{\mb{\Lambda}_1}> \norm{\mb{\Lambda}_2}\) since \(\mb{\Lambda}\cdot\mb{\Lambda}\) is a negative scalar.
I will call \(\mb{\Lambda}_1\) the larger axis and  \(\mb{\Lambda}_2\) the smaller axis.

If \((\mb{\Lambda}\cdot\mb{\Lambda})^2=(\mb{\Lambda}\vee\mb{\Lambda})^2\) 
or equivalently \(\mb{\Lambda}\cdot\mb{\Lambda}=\pm\mb{\Lambda}\vee\mb{\Lambda}\),
then 
\begin{equation}
\mb{\Lambda}_1^2=\mb{\Lambda}_2^2=\tfrac{1}{2}(\mb{\Lambda}\cdot\mb{\Lambda}).
\end{equation}
These two cases exhaust all possibilities, since
\((\mb{\Lambda}\cdot\mb{\Lambda})^2\ge(\mb{\Lambda}\vee\mb{\Lambda})^2\) for any bivector \(\mb{\Lambda}\) in \El{3}.
Hence, if \(\mb{\Lambda}_1^2\) and \(\mb{\Lambda}_2^2\) are distinct, the axes  of \(\mb{\Lambda}\) are given by
\begin{equation}
\mb{\Lambda}_1=\mb{\Lambda}/(1+\tfrac{1}{2}(\mb{\Lambda}\wedge\mb{\Lambda})/\mb{\Lambda}_1^2)\quad
\text{and}\quad
\mb{\Lambda}_2=\mb{\Lambda}/(1+\tfrac{1}{2}(\mb{\Lambda}\wedge\mb{\Lambda})/\mb{\Lambda}_2^2).
\end{equation}
Note that \(\mb{\Lambda}_2/\norm{\mb{\Lambda}_2}=\mb{\Lambda}_1\I/\norm{\mb{\Lambda}_1}\)
if \(\mb{\Lambda}\vee\mb{\Lambda}<0\) and 
\(\mb{\Lambda}_2/\norm{\mb{\Lambda}_2}=-\mb{\Lambda}_1\I/\norm{\mb{\Lambda}_1}\)
if \(\mb{\Lambda}\vee\mb{\Lambda}>0\),
so the axes are at the distance \(\tfrac{\pi}{2}\) from each other.

If \(\mb{\Lambda}_1^2=\mb{\Lambda}_2^2\), 
I get \(\mb{\Lambda}_{1,2}(1+\mb{\Lambda}\wedge\mb{\Lambda}/\mb{\Lambda}\cdot\mb{\Lambda})=\mb{\Lambda}\)
but the multivector \(1+\mb{\Lambda}\wedge\mb{\Lambda}/\mb{\Lambda}\cdot\mb{\Lambda}=1\pm\I\) is not invertible.
In this case,
the decomposition of the non-simple \(\mb{\Lambda}\) into complementary lines is not unique.
If \(\mb{\Lambda}\cdot\mb{\Lambda}=\mb{\Lambda}\vee\mb{\Lambda}\), then 
\begin{equation}
\mb{\Lambda}=\lambda_1\e_{10}+\lambda_2\e_{20}+\lambda_3\e_{30}+\lambda_1\e_{23}+\lambda_2\e_{31}+\lambda_3\e_{12}
\end{equation}
for some \(\lambda_1, \lambda_2, \lambda_3\in\R{}\) and \(\mb{\Lambda}\) is a positive Clifford bivector.
One possible decomposition is given by 
\(\mb{\Lambda}_1=\lambda_1\e_{23}+\lambda_2\e_{31}+\lambda_3\e_{12}\) and 
\(\mb{\Lambda}_2=\mb{\Lambda}_1\I=\lambda_1\e_{10}+\lambda_2\e_{20}+\lambda_3\e_{30}\).
The bivector \(\mb{\Lambda}\) can also be decomposed as 
\(\mb{\Lambda}=\mb{\Lambda}_1^\p+\mb{\Lambda}_1^\p\I\), where
\(\mb{\Lambda}_1^\p\) is any positive  Clifford parallel of \(\mb{\Lambda}_1\)
(one needs to ensure that the weight of \(\mb{\Lambda}_1^\p\) matches that of \(\mb{\Lambda}_1\)).
On the other hand,
if \(\mb{\Lambda}\cdot\mb{\Lambda}=-\mb{\Lambda}\vee\mb{\Lambda}\), then 
\begin{equation}
\mb{\Lambda}=-\lambda_1\e_{10}-\lambda_2\e_{20}-\lambda_3\e_{30}+\lambda_1\e_{23}+\lambda_2\e_{31}+\lambda_3\e_{12}
\end{equation}
for some \(\lambda_1, \lambda_2, \lambda_3\in\R{}\) and \(\mb{\Lambda}\) is a negative Clifford bivector.
A decomposition is given by 
\(\mb{\Lambda}_1=\lambda_1\e_{23}+\lambda_2\e_{31}+\lambda_3\e_{12}\) and 
\(\mb{\Lambda}_2=-\mb{\Lambda}_1\I=-\lambda_1\e_{10}-\lambda_2\e_{20}-\lambda_3\e_{30}\).
Other decompositions are given by \(\mb{\Lambda}=\mb{\Lambda}_1^\n-\mb{\Lambda}_1^\n\I\),
where \(\mb{\Lambda}_1^\n\) is any negative Clifford parallel of \(\mb{\Lambda}_1\).

\begin{figure}[t]
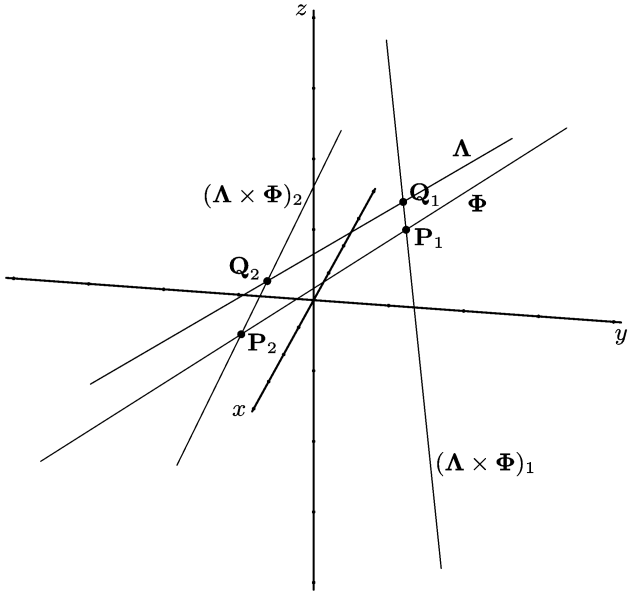
\hspace{-1cm}
\begin{asy}
import Figure3D;
Figure f = Figure();
metric=Metric(Elliptic);

MV L = Line(-3/2,1,-1/2,-1,-5/2,-2); //Line(-3/2,1,-1,-1,-5/2,-1);
MV K = Line(1,5/3,-2,-1,3,2); //Line(1,4/3,-2,0,3,2);

f.line(L, "$\boldsymbol{\Lambda}$", position=0.1,align=(0,-0.5,1),draw_orientation=false);
f.line(K, "$\boldsymbol{\Phi}$", position=0.8,align=(0,1,-0.5),draw_orientation=false);

MV F = cross(L,K);
MV F1=axis1(F);
MV F2=axis2(F);

f.line(F1,"$(\boldsymbol{\Lambda}\times\boldsymbol{\Phi})_1$",draw_orientation=false,size=4);
f.line(F2,"$(\boldsymbol{\Lambda}\times\boldsymbol{\Phi})_2$",position=0.8,align=(0,-1,0),draw_orientation=false,size=4.5);

MV Q1 = wedge(F1,join(L,O)); Q1/=norm(Q1);
MV P1 = wedge(F1,join(K,O)); P1/=norm(P1);

f.point(Q1,"$\textbf{Q}_1$",align=(0,1,0.5),draw_orientation=false,draw_helper_lines=false);
f.point(P1,"$\textbf{P}_1$",align=(0,1.25,-0.5),draw_orientation=false,draw_helper_lines=false);
asin(norm(join(P1,Q1)));

MV Q2 = wedge(F2,join(L,O)); Q2/=norm(Q2);
MV P2 = wedge(F2,join(K,O)); P2/=norm(P2);

f.point(Q2,"$\textbf{Q}_2$",align=(0,-1.25,0.75),draw_orientation=false,draw_helper_lines=false);
f.point(P2,"$\textbf{P}_2$",align=(0,1,-0.5),draw_orientation=false,draw_helper_lines=false);
asin(norm(join(P2,Q2)));

\end{asy}
\caption{Commutator in \El{3}}
\label{commutator in El3}
\end{figure}

The commutator \(\mb{\Lambda}\times\mb{\Phi}\) of two intersecting lines \(\mb{\Lambda}\) and \(\mb{\Phi}\) is simple.
It is a line which is perpendicular to both \(\mb{\Lambda}\) and \(\mb{\Phi}\)
and passes through \(\mb{\Lambda}\) and \(\mb{\Phi}\)  at the point where they intersect.
On the other hand, if \(\mb{\Lambda}\) and \(\mb{\Phi}\) are Clifford-parallel,
then  \(\mb{\Lambda}\times\mb{\Phi}\) is generally not simple but the decomposition into axes is not unique.
In fact, for any point on \(\mb{\Lambda}\) (or on \(\mb{\Phi}\)), it is possible
to find an axis which passes through the point. 
The separation \(r\) between  \(\mb{\Lambda}\) and \(\mb{\Phi}\)  measured along the axes,
i.e.\ the distance between the Clifford-parallel lines \(\mb{\Lambda}\) and \(\mb{\Phi}\), is constant
and satisfies \(\sin^2 r=|\mb{\Lambda}\vee\mb{\Phi}|\) and \(\cos^2 r=|\mb{\Lambda}\cdot\mb{\Phi}|\)
where \(\mb{\Lambda}\) and  \(\mb{\Phi}\) are assumed to be normalised.

If \(\mb{\Lambda}\) and \(\mb{\Phi}\) do not intersect and are not Clifford-parallel,
then \(\mb{\Lambda}\times\mb{\Phi}\) is non-simple and can be decomposed uniquely into two axes,
\((\mb{\Lambda}\times\mb{\Phi})_1\) and \((\mb{\Lambda}\times\mb{\Phi})_2\).
The axes are perpendicular to both  \(\mb{\Lambda}\) and \(\mb{\Phi}\)
and pass through both lines.
Since the axes are perpendicular to both \(\mb{\Lambda}\) and \(\mb{\Phi}\),
the separation between \(\mb{\Lambda}\) and \(\mb{\Phi}\) measured along each axis is at a local minimum.
If \(r_1\) and \(r_2\) denote these local minimal separations, then \(r_1<r_2\) and 
the distance \(r\) between the lines \(\mb{\Lambda}\) and \(\mb{\Phi}\)
is given by \(r=r_1\), where \(r_1\) corresponds to the separation along the larger axis \((\mb{\Lambda}\times\mb{\Phi})_1\)
and \(r_2\) corresponds to the separation along the smaller axis \((\mb{\Lambda}\times\mb{\Phi})_1\).

An example is shown in Figure~\ref{commutator in El3}, where 
\(\mb{\Lambda}=-\tfrac{3}{2}\e_{10}+\e_{20}-\tfrac{1}{2}\e_{30}-\e_{23}-\tfrac{5}{2}\e_{31}-2\e_{12}\) 
and \(\mb{\Phi}=\e_{10}+\tfrac{5}{3}\e_{20}-2\e_{30}-\e_{23}+3\e_{31}+2\e_{12}\).
I get 
\(\sin r=\sin r_1=\norm{\tb{P}_1\vee\tb{Q}_1}=\scriptstyle \sqrt{(22-5\sqrt{17})/59}\) 
for the normalised points \(\tb{P}_1\) and \(\tb{Q}_1\),
where the axis \((\mb{\Lambda}\times\mb{\Phi})_1\) intersects  the lines \(\mb{\Lambda}\) and \(\mb{\Phi}\),
and 
\(\sin r_2=\norm{\tb{P}_2\vee\tb{Q}_2}=\scriptstyle \sqrt{(22+5\sqrt{17})/59}\) 
for the normalised points \(\tb{P}_2\) and \(\tb{Q}_2\),
where \((\mb{\Lambda}\times\mb{\Phi})_2\) intersects  \(\mb{\Lambda}\) and \(\mb{\Phi}\).

The local minimal separations \(r_1\) and \(r_2\) satisfy
\(\cos r_1\cos r_2 = |\mb{\Lambda}\cdot\mb{\Phi}|\) and
\(\sin r_1\sin r_2 = |\mb{\Lambda}\vee\mb{\Phi}|\),
provided that the lines are normalised.
The rotation around the axis \((\mb{\Lambda}\times\mb{\Phi})_2\) causes the point \(\tb{Q}_1\)
to moves along \((\mb{\Lambda}\times\mb{\Phi})_1\).
The distance \(r=r_1\) can be thought of as the angle that the point \(\tb{Q}_1\) needs to be rotated by
in order to bring it to the point \(\tb{P}_1\)
(note that \(\tb{Q}_1\) must be rotated in the correct direction to get \(r\);
if \(\tb{Q}_1\) is rotated in the opposite direction one gets \(\pi-r\) instead).
Similarly, \(r_2\) can be thought of as the angle required to bring \(\tb{Q}_2\) to \(\tb{P}_2\)
when the point \(\tb{Q}_2\) is rotated around the axis \((\mb{\Lambda}\times\mb{\Phi})_1\).
It is consistent with the angle between two planes, both of which pass through
the axis \((\mb{\Lambda}\times\mb{\Phi})_1\)
and one of which contains \(\mb{\Lambda}\) while the other contains \(\mb{\Phi}\).

By analogy with the Euclidean case, one can refer to \(r_2\) as the angle between the lines 
\(\mb{\Lambda}\) and \(\mb{\Phi}\). 
Indeed, in the Euclidean case, the distance between two lines can be obtained by the rotation
around the axis at infinity, i.e.\ translation, while the angle can be obtained by the rotation
around the finite axis.
In the elliptic case, both axes are finite and it is natural to associate the smaller separation with the distance
and the larger separation with the angle between the lines.
Since \(r_2\in[0,\tfrac{\pi}{2}]\), it does not take the orientation of the lines into account.
It is more convenient to define the angle \(\alpha\in[0,\pi]\) between the normalised lines \(\mb{\Lambda}\) and \(\mb{\Phi}\)
by \(\mb{\Lambda}\cdot\mb{\Phi}=-\cos r\cos\alpha\), where \(r\) is the distance between the lines,
in order to take the orientation of the lines into account
(\(r_2\) may be called the unoriented angle 
and \(r_2=\alpha\) if \(\alpha\le\tfrac{\pi}{2}\) or \(r_2=\pi-\alpha\) if \(\alpha>\tfrac{\pi}{2}\)).
Hence, in general the geometric product of two normalised lines \(\mb{\Lambda}\) and \(\mb{\Phi}\) is given by
\begin{equation}
\mb{\Lambda}\mb{\Phi}=-\cos r\cos\alpha+\mb{\Lambda}\times\mb{\Phi}\pm\I\sin r\sin\alpha,
\end{equation}
where the sign in front of \(\I\) depends on the relative orientation of the lines.
The distance \(r\in[0,\tfrac{\pi}{2}]\)  between normalised lines \(\mb{\Lambda}\) and \(\mb{\Phi}\)
can thus be determined from
\begin{equation}
2\sin^2r =1+v^2-u^2-\sqrt{(1+v^2-u^2)^2-4v^2},
\end{equation}
where  \(u=\mb{\Lambda}\cdot\mb{\Phi}\) and \(v=\mb{\Lambda}\vee\mb{\Phi}\).
Once the distance is found, the  angle \(\alpha\in[0,\pi]\)  is determined from \(\cos\alpha=-u/\cos r\).
Since both distance and angular measures in \El{3} are elliptic, there is not much difference between angles and distances
as evidenced by the above discussion.

\begin{figure}[t]\hspace{-1cm}
\begin{subfloatenv}{\El{3}}
\begin{asy}
import Figure3D;
Figure f = Figure();
metric=Metric(Elliptic);

var a = Plane(1,1/2,-3/2,1);
f.plane(a, "$\textbf{a}$", align=(1,-3,-1),draw_orientation=false);

var P = Point(1,-1/2,0,1.5); //Point(1,-1,1,1.5);
f.point(P, "$\textbf{P}$", draw_orientation=false, align=(0,0.5,1));

f.point(a*I, "$\textbf{a}\textbf{I}$", draw_orientation=false, align=(0,1,-1));
f.line(dot(a,P), "$\textbf{a}\cdot\textbf{P}$", position=0.02, draw_orientation=false, align=(0,-1,1.5));

f.point(dot(P,a)/a, "$(\textbf{P}\cdot\textbf{a})\textbf{a}^{-1}$", draw_orientation=false, align=(0,-0.5,0.5));
f.point(wedge(P,a)/a, "$(\textbf{P}\wedge\textbf{a})\textbf{a}^{-1}$", draw_orientation=false, align=(0,-0.5,0.5));

\end{asy}
\end{subfloatenv}\hfill%
\begin{subfloatenv}{\El{3}}
\begin{asy}
import Figure3D;
Figure f = Figure();
metric=Metric(Elliptic);

var a = Plane(1,1/2,-3/2,1);
f.plane(a, "$\textbf{a}$", align=(1,-3,-1),draw_orientation=false);
f.point(a*I, "$\textbf{a}\textbf{I}$", draw_orientation=false, align=(0,1,-1));

MV K = join(Point(1,-1,1,1),Point(1,2,-2,-1));
f.line(K,"$\boldsymbol{\Lambda}$",draw_orientation=false,align=(0,-0.5,1));
f.line(dot(K,a)/a,"$(\boldsymbol{\Lambda}\cdot\textbf{a})\textbf{a}^{-1}$",align=(0,-1,1),position=0.2,draw_orientation=false);
f.line(wedge(K,a)/a,"$(\boldsymbol{\Lambda}\wedge\textbf{a})\textbf{a}^{-1}$",align=(0,-0.5,1),position=0.25,draw_orientation=false);
f.line(a*K/a,"$\textbf{a}\boldsymbol{\Lambda}\textbf{a}^{-1}$",align=(0,1,1),position=0.95,draw_orientation=false,size=3);

\end{asy}
\end{subfloatenv}
\caption{Projection, rejection and reflection in \El{3} (1)}
\label{projection and rejection in El3}
\end{figure}

\begin{figure}[t]
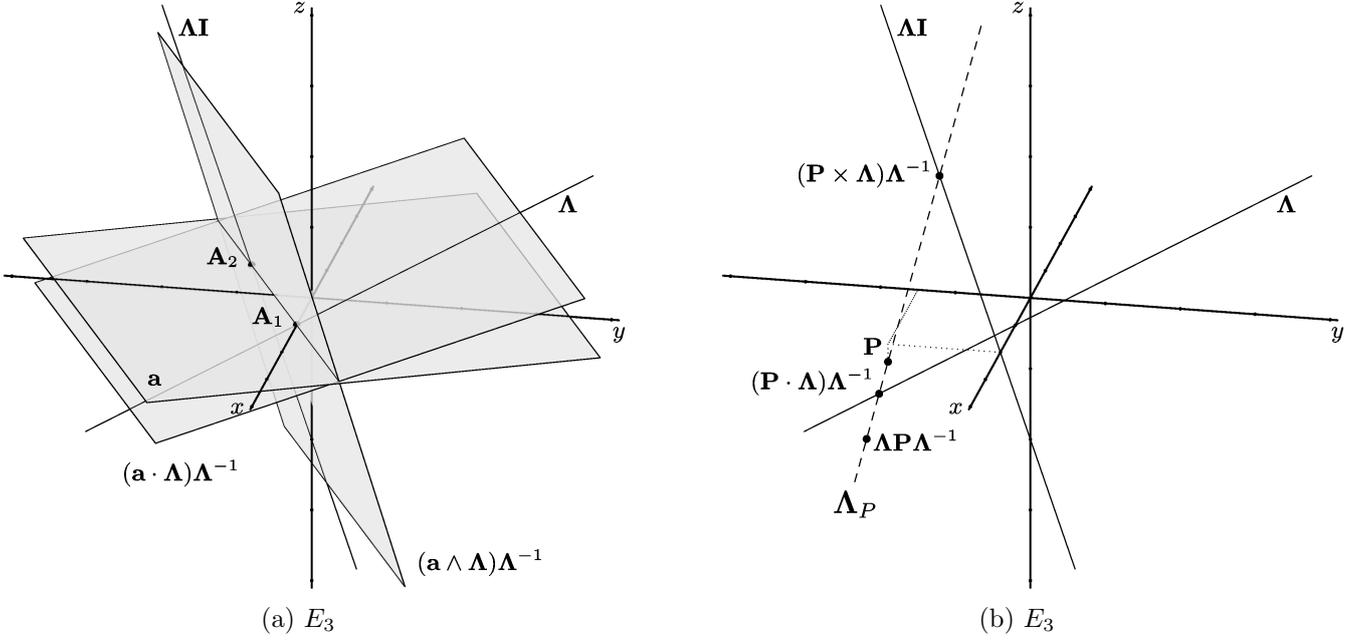
\hspace{-1cm}
\begin{subfloatenv}{\El{3}}
\begin{asy}
import Figure3D;
Figure f = Figure();
metric=Metric(Elliptic);

var P = Point(1,1,0,0);
var Q = Point(1,0,1,1/3);
var L = join(P,Q);
write(L);

f.line(L, "$\boldsymbol{\Lambda}$", position=0.05,align=(0,0,-1),draw_orientation=false);
f.line(L*I, "$\boldsymbol{\Lambda}\textbf{I}$", position=0.95,align=(0,1,0.5),draw_orientation=false);

var a = e_3;
f.point(wedge(a,L),"$\textbf{A}_1$",align=(1,-2,1),draw_orientation=false,draw_helper_lines=false);
f.point(wedge(a,L*I),"$\textbf{A}_2$",align=(1,-2,1),draw_orientation=false,draw_helper_lines=false);

MV K = wedge(wedge(a,L)/L,dot(a,L)/L);
triple u = (K.p23,K.p31,K.p12);
f.plane(a, "$\textbf{a}$",u=u,align=(1,1,3),draw_orientation=false);
f.plane(wedge(a,L)/L, "$(\textbf{a}\wedge\boldsymbol{\Lambda})\boldsymbol{\Lambda}^{-1}$", u=u, size=3,align=(0,2,2),label_position=3,draw_orientation=false);
f.plane(dot(a,L)/L, "$(\textbf{a}\cdot\boldsymbol{\Lambda})\boldsymbol{\Lambda}^{-1}$", u=u, size=3,align=(1,1,-2),draw_orientation=false);

//f.plane(L*a/L, "$\boldsymbol{\Lambda}\textbf{a}\boldsymbol{\Lambda}^{-1}$", u=u, size=3,align=(1,1,-2),draw_orientation=false);
\end{asy}
\end{subfloatenv}\hfill%
\begin{subfloatenv}{\El{3}}
\begin{asy}
import Figure3D;
Figure f = Figure();
metric=Metric(Elliptic);

var LL1 = Point(1,1,0,0);
var LL2 = Point(1,0,1,1/3);
var L = join(LL1,LL2);
write(L);

f.line(L, "$\boldsymbol{\Lambda}$", position=0.05,align=(0,0,-1),draw_orientation=false);
f.line(L*I, "$\boldsymbol{\Lambda}\textbf{I}$", position=0.95,align=(0,1,0.5),draw_orientation=false);

//MV C = toline(L).centre();
//f.point(C," ",draw_orientation=false,draw_helper_lines=true);

//MV CI = toline(L*I).centre();
//f.point(CI," ",draw_orientation=false,draw_helper_lines=true);

MV P = Point(1,2,-1.5,-1/4); //Point(1,5/3,-1.5,0);// Point(1,-0,-1.5,1);
f.point(P, "$\textbf{P}$",align=(0,-1,1),draw_orientation=false,draw_helper_lines=true);
f.point(cross(P,L)/L, "$(\textbf{P}\times\boldsymbol{\Lambda})\boldsymbol{\Lambda}^{-1}$",align=(0,-1.5,0),draw_orientation=false,draw_helper_lines=false);
f.point(dot(P,L)/L, "$(\textbf{P}\cdot\boldsymbol{\Lambda})\boldsymbol{\Lambda}^{-1}$",align=(1,-1,1),draw_orientation=false,draw_helper_lines=false);

 f.point(L*P/L, "$\boldsymbol{\Lambda}\textbf{P}\boldsymbol{\Lambda}^{-1}$",align=(0,1,0),draw_orientation=false,draw_helper_lines=false);

MV Lp = wedge(dot(L,P),join(L,P));
//f.line(Lp,pen=dashed,draw_orientation=false);
f.line(Lp,"$\boldsymbol{\Lambda}_P$",position=0,align=(0,0,-1),pen=dashed,draw_orientation=false);

//MV P1 = wedge(P*I,L);
//MV P2 = wedge(P*I,L*I);
//f.point(P1,"$\textbf{P}_1$",align=(0,-0.25,1.25),draw_orientation=false,draw_helper_lines=false);
//f.point(P2,"$\textbf{P}_2$",align=(0,-0.5,-1),draw_orientation=false,draw_helper_lines=false);
//MV K = join(P1,P2);
//f.line(K*I,pen=dashed,draw_orientation=false);
//f.line(K,pen=dashed,draw_orientation=false,size=5) ;

join(wedge(P*I,L),wedge(P*I,L*I))*I -0.5*join(L*P*L,P);

\end{asy}
\end{subfloatenv}
\caption{Projection, rejection and reflection in \El{3} (2)}
\label{projection and rejection in El3 lines}
\end{figure}

The usual formulas for projection and rejection apply in \El{3}.
The projection of a geometric object dually represented by a blade 
\(B_l\) on a plane \(\tb{a}\) is defined by \((B_l\cdot\tb{a})\tb{a}^{-1}\)
and the rejection is defined by \((B_l\wedge\tb{a})\tb{a}^{-1}\).
For instance the rejection of a line \(\mb{\Lambda}\) by \(\tb{a}\) is given by \((\mb{\Lambda}\wedge\tb{a})\tb{a}^{-1}\).
The rejections pass through the polar point \(\tb{a}\I\) of the plane \(\tb{a}\),
and the projections lie on \(\tb{a}\).
The rejection by \(\tb{a}\) is zero if \(B_l\) lies on the plane \(\tb{a}\), and the projection on \(\tb{a}\)
is zero if \(B_l\) passes through 
the polar point \(\tb{a}\I\) of \(\tb{a}\), in which case \(B_l\) is perpendicular to \(\tb{a}\).
An example of the projection on and rejection by the plane \(\tb{a}=\e_0+\tfrac{1}{2}\e_1-\tfrac{3}{2}\e_2+\e_3\)
is shown in Figure~\ref{projection and rejection in El3}
for the point \(\tb{P}=\e_{123}-\tfrac{1}{2}\e_{320}+\tfrac{3}{2}\e_{210}\) and 
the line  \(\mb{\Lambda}=\e_{10}+\e_{20}-3\e_{23}+3\e_{31}+2\e_{12}\).

The projection of \(B_l\) on a point \(\tb{P}\) is defined by \((B_l\cdot\tb{P})\tb{P}^{-1}\)
and the rejection is defined by \((B_l\wedge\tb{P})\tb{P}^{-1}\) for planes \((l=1)\) and by 
 \((B_l\times\tb{P})\tb{P}^{-1}\) for lines and points \((l=2,3)\), in accord with the decomposition of the relevant geometric products.
For instance, the rejection of a line \(\mb{\Lambda}\) by \(\tb{P}\) is given by \((\mb{\Lambda}\times\tb{P})\tb{P}^{-1}\).
The rejections lie on the plane \(\tb{P}\I\), whose polar point is \(\tb{P}\), 
and the projections pass through the point \(\tb{P}\).
The projection of \(B_l\) on \(\tb{P}\) is zero if \(B_l\) lies on \(\tb{P}\I\),
and the rejection of \(B_l\) by \(\tb{P}\) is zero if \(B_l\) passes through \(\tb{P}\).
The configurations are in some sense dual to those arising in projection on and rejection by a plane.

The projection of \(B_l\) on a line \(\mb{\Lambda}\) is defined by \((B_l\cdot\mb{\Lambda})\mb{\Lambda}^{-1}\) 
for points and planes \((l=1,3)\).
The rejection of a plane \(\tb{a}\) by \(\mb{\Lambda}\) is defined by \((\tb{a}\wedge\mb{\Lambda})\mb{\Lambda}^{-1}\),
and the rejection of a point \(\tb{P}\) by \(\mb{\Lambda}\) is defined by \((\tb{P}\times\mb{\Lambda})\mb{\Lambda}^{-1}\).
An example of the projection on and rejection by the line 
\(\mb{\Lambda}=-\tfrac{1}{3}\e_{20}+\e_{30}+\e_{23}-\e_{31}-\tfrac{1}{3}\e_{12}\)
 is shown in Figure~\ref{projection and rejection in El3 lines}
for the plane \(\tb{a}=\e_{3}\) and the point
\(\tb{P}=\e_{123}+2\e_{320}-\tfrac{3}{2}\e_{130}-\tfrac{1}{4}\e_{210}\).
In Figure~\ref{projection and rejection in El3 lines}(a), 
the points \(\tb{A}_1\) and \(\tb{A}_2\) are located at the intersect of the plane \(\tb{a}\)
with the lines \(\mb{\Lambda}\) and \(\mb{\Lambda}\I\), respectively.
The projection \((\tb{a}\cdot\mb{\Lambda})\mb{\Lambda}^{-1}\) contains \(\mb{\Lambda}\)
and the rejection \((\tb{a}\wedge\mb{\Lambda})\mb{\Lambda}^{-1}\) contains \(\mb{\Lambda}\I\).
The projection \((\tb{P}\cdot\mb{\Lambda})\mb{\Lambda}^{-1}\)
and rejection \((\tb{P}\times\mb{\Lambda})\mb{\Lambda}^{-1}\) of the point \(\tb{P}\)
by the line \(\mb{\Lambda}\) both lie on 
the line \(\mb{\Lambda}_P=(\mb{\Lambda}\cdot\tb{P})\wedge(\mb{\Lambda}\vee\tb{P})\), which
passes through \(\tb{P}\) and is perpendicular to \(\mb{\Lambda}\);
\(\mb{\Lambda}_P\) is also perpendicular to \(\mb{\Lambda}\I\).
Note that \(\mb{\Lambda}_P\) can also be written as 
\(\mb{\Lambda}_P=\tfrac{1}{2}\tb{P}\vee(\mb{\Lambda}\tb{P}\mb{\Lambda}^{-1})\)
if \(\mb{\Lambda}\) is normalised.

The geometric product  \(\mb{\Phi}\mb{\Lambda}\) of two lines 
consists of three terms which can be split into two components
for the projection and rejection in two different ways.
In general, the commutator \(\mb{\Phi}\times\mb{\Lambda}\) is not simple, so it needs to be decomposed into two axes.
Then the larger axis \((\mb{\Phi}\times\mb{\Lambda})_1\) can be combined with \(\mb{\Phi}\cdot\mb{\Lambda}\)
and the smaller axis \((\mb{\Phi}\times\mb{\Lambda})_2\) with \(\mb{\Phi}\wedge\mb{\Lambda}\),
which gives the first kind of projection and rejection:
\begin{equation}
\begin{split}
&\mathsf{proj_1}(\mb{\Phi};\mb{\Lambda})=
(\mb{\Phi}\cdot\mb{\Lambda}+(\mb{\Phi}\times\mb{\Lambda})_1)\mb{\Lambda}^{-1},\\
&\mathsf{rej_1}(\mb{\Phi};\mb{\Lambda})=
((\mb{\Phi}\times\mb{\Lambda})_2+\mb{\Phi}\wedge\mb{\Lambda})\mb{\Lambda}^{-1}.
\end{split}
\end{equation}
If \((\mb{\Phi}\times\mb{\Lambda})_1\) is combined with \(\mb{\Phi}\wedge\mb{\Lambda}\)
and \((\mb{\Phi}\times\mb{\Lambda})_2\) with \(\mb{\Phi}\cdot\mb{\Lambda}\), 
one gets the second kind of projection and rejection:
\begin{equation}
\begin{split}
&\mathsf{proj_2}(\mb{\Phi};\mb{\Lambda})=
(\mb{\Phi}\cdot\mb{\Lambda}+(\mb{\Phi}\times\mb{\Lambda})_2)\mb{\Lambda}^{-1},\\
&\mathsf{rej_2}(\mb{\Phi};\mb{\Lambda})=
((\mb{\Phi}\times\mb{\Lambda})_1+\mb{\Phi}\wedge\mb{\Lambda})\mb{\Lambda}^{-1}.
\end{split}
\end{equation}
These projections and rejections share some of their basic properties with those in Euclidean space.
For instance, 
\(\mathsf{proj_1}(\mb{\Phi};\mb{\Lambda})\)
is a line which passes through \(\mb{\Lambda}\) and is at the same angle to \(\mb{\Lambda}\) as \(\mb{\Phi}\), and
\(\mathsf{rej_2}(\mb{\Phi};\mb{\Lambda})\)
is a line which is at the same distance from \(\mb{\Lambda}\) as \(\mb{\Phi}\) and
is perpendicular to \(\mb{\Lambda}\).

%
%
%
%
%
%
%
%

The top-down reflection of \(B_l\) in \(A_k\) is given by \((-1)^{kl}A_kB_lA_k^{-1}\)
and the bottom-up reflection is given by  \((-1)^{k(l-1)}A_kB_lA_k^{-1}\).
The top-down and bottom-up reflections in lines are identical since \(k=2\) for lines.
For instance, the reflection of a point \(\tb{P}\) in a line \(\mb{\Lambda}\) is given by 
\(\mb{\Lambda}\tb{P}\mb{\Lambda}^{-1}\),
and the top-down reflection of \(\mb{\Lambda}\) in a plane \(\tb{a}\) is given by  \(\tb{a}\mb{\Lambda}\tb{a}^{-1}\).
Some examples are shown in Figures~\ref{projection and rejection in El3} and \ref{projection and rejection in El3 lines}.
Note that the points \(\tb{P}\) 
and  \(\mb{\Lambda}\tb{P}\mb{\Lambda}^{-1}\) are at the same distance from \(\mb{\Lambda}\) 
(see Figure~\ref{projection and rejection in El3 lines}(b)),
and
the line \(\tb{a}\mb{\Lambda}\tb{a}^{-1}\) is at the same angle to the plane \(\tb{a}\) as \(\mb{\Lambda}\)
(see Figure~\ref{projection and rejection in El3}(b)).

Spinors in \El{3} are defined in the usual way as multivectors that can be written as the product of an even number
of planes which square to unity.
They are even and each spinor satisfies \(S\reverse{S}=1\).
Any spinor in \El{3} can also be written as \(S=e^A\) for some bivector \(A\).
Spinors form a Lie group whose Lie algebra consists of bivectors with the commutator used as the product in the algebra.

\begin{figure}[t]
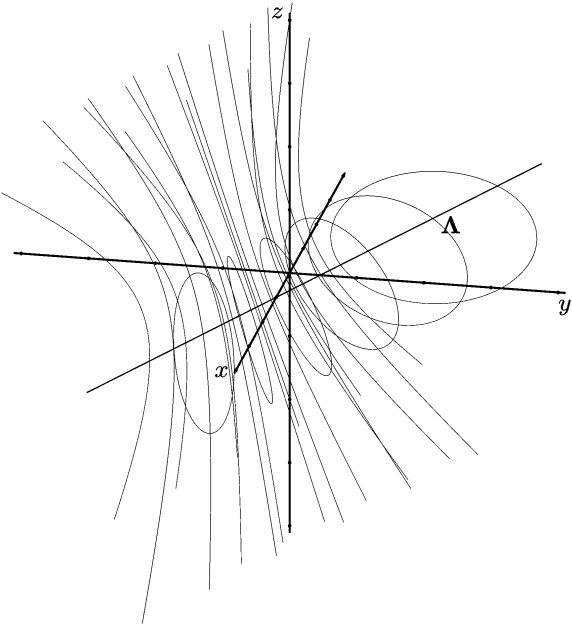
\hspace{-1cm}
\begin{asy}
import Figure3D;
Figure f = Figure();
metric=Metric(Elliptic);

var P = Point(1,1,0,0);
var Q = Point(1,0,1,1/3);
var L = join(P,Q);
L/=norm(L);
write(L);

f.line(L, "$\boldsymbol{\Lambda}$", align=(0,0,-1),draw_orientation=false);

MV C = toline(L).centre();
triple n = (L.p23,L.p31,L.p12);
n = normalise(n);

MV K = Line(0,0,0,L.p23,L.p31,L.p12); K/=norm(K);
MV M = 1.2*K*L/K - 0.2*K; M/=norm(M);

MV[] points = {};
MV[] pointsK = {};
for(int i=-2; i<=3; ++i) { triple p = totriple(C)+0.6*i*n; MV P = Point(1,p.x,p.y,p.z); points.push(P); pointsK.push(wedge(K,join(P,L*I))); }
MV[] pointsM = {};
for(int i=-2; i<=3; ++i) { triple p = totriple(C)+1.2*i*n; MV P = Point(1,p.x,p.y,p.z); points.push(P); pointsM.push(wedge(M,join(P,L*I))); }

//for(MV p: pointsK) { f.point(p," ",draw_orientation=false,draw_helper_lines=false); }
//for(MV p: pointsM) { f.point(p," ",draw_orientation=false,draw_helper_lines=false); }

//triple u = ((L*I).p23,(L*I).p31,(L*I).p12);
//for(MV P: points) { f.plane(join(P,L*I),u=u,draw_orientation=false); }

for(MV P: pointsK) {
  triple f(real t) { MV S = cos(1/2*t)-L*sin(1/2*t); return totriple(S*P/S); };
  real x(real t) { return f(t).x; };
  real y(real t) { return f(t).y; };
  real z(real t) { return f(t).z; };
  path3 c = graph(x, y, z, -pi, pi, operator ..);
  draw(c, currentpen+0.25);
}

for(MV P: pointsM) {
  triple f(real t) { MV S = cos(1/2*t)-L*sin(1/2*t); return totriple(S*P/S); };
  real x(real t) { return f(t).x; };
  real y(real t) { return f(t).y; };
  real z(real t) { return f(t).z; };
  path3 c = graph(x, y, z, -pi/2+0.17, pi/2-0.17, operator ..);
  draw(c, currentpen+0.25);
}

for(MV P: pointsM) {
  triple f(real t) { MV S = cos(1/2*t)-L*sin(1/2*t); return totriple(S*P/S); };
  real x(real t) { return f(t).x; };
  real y(real t) { return f(t).y; };
  real z(real t) { return f(t).z; };
  path3 c = graph(x, y, z, pi/2+0.78, 3*pi/2-0.78, operator ..);
  draw(c, currentpen+0.25);
}

for(MV P: pointsM) {
  triple f(real t) { MV S = cos(1/2*t)-L*sin(1/2*t); return totriple(S*P/S); };
  real x(real t) { return f(t).x; };
  real y(real t) { return f(t).y; };
  real z(real t) { return f(t).z; };
  path3 c = graph(x, y, z, -pi, -pi/2-0.6, operator ..);
  draw(c, currentpen+0.25);
}

MV C = toline(L*I).centre();
MV[] pointsLI = {};
for(int i=-2; i<=3; ++i) { triple p = totriple(C)+0.6*i*n; MV P = Point(1,p.x,p.y,p.z); pointsLI.push(P); }

for(MV P: pointsLI) {
  triple f(real t) { MV S = cos(1/2*t)-L*sin(1/2*t); return totriple(S*P/S); };
  real x(real t) { return f(t).x; };
  real y(real t) { return f(t).y; };
  real z(real t) { return f(t).z; };
  path3 c = graph(x, y, z, -1, 1, operator ..);
  draw(c, currentpen+0.25);
}

\end{asy}
\caption{Rotation in \El{3}}
\label{rotation in El3}
\end{figure}

Any proper motion in \El{3} can be obtained as an action of a spinor.
One encounters the following proper motions.
The action \(SA_kS^{-1}\) of a spinor 
\begin{equation}
S=e^{-\tfrac{1}{2}(\alpha+\beta\I)\mb{\Lambda}},
\end{equation}
where \(\alpha,\beta\in\R{}\), \(\alpha\ne\beta\), and \(\mb{\Lambda}\) is a normalised line,
on the geometric object dually represented by a blade \(A_k\)
yields a double rotation of \(A_k\) around \(\mb{\Lambda}\) by the angle \(\alpha\)
and around \(\mb{\Lambda}\I\) by \(\beta\).
Under the double rotation, the axes  \(\mb{\Lambda}\) and  \(\mb{\Lambda}\I\) of the double rotation are invariant,
which implies that any point on \(\mb{\Lambda}\) or  \(\mb{\Lambda}\I\) stays on the line.
Double rotation in \El{3} corresponds to a proper motion in \E{3} which consists of a simple rotation 
around an axis and a translation along the same axis.
If \(\beta=0\), a simple rotation around \(\mb{\Lambda}\) is obtained.
Some trajectories generated by the rotation around the line \(\mb{\Lambda}\) are shown in Figure~\ref{rotation in El3}.

If one substitutes \(\alpha=\beta\) in the above, the spinor can be written as
\begin{equation}
S=e^{-\tfrac{1}{2}\beta\mb{\Xi}^+},
\label{Clifford translation 1}
\end{equation}
where \(\mb{\Xi}^+=(\I+1)\mb{\Lambda}\) is a positive Clifford bivector,
and the resulting proper motion occurs along the Clifford parallels associated with \(\mb{\Xi}^+\).
If, on the other hand, \(\alpha=-\beta\), then
\begin{equation}
S=e^{-\tfrac{1}{2}\beta\mb{\Xi}^-},
\label{Clifford translation 2}
\end{equation}
where \(\mb{\Xi}^-=(\I-1)\mb{\Lambda}\) is a negative Clifford bivector,
and the proper motion is along the Clifford parallels of \(\mb{\Xi}^-\).
The proper motions which result from the action \(SA_kS^{-1}\) of \(S\) defined by 
(\ref{Clifford translation 1}) and (\ref{Clifford translation 2})
are called  Clifford translations.
The Clifford translation of a point \(\tb{P}\) generated by  a Clifford bivector \(\mb{\Xi}\) can also be written as 
\begin{equation}
\tb{P}\mapsto \tb{P}\cos\beta +\tb{P}\times\mb{\Xi}\sin\beta
\end{equation}
by expressing \(\mb{\Xi}\) in terms of the Clifford parallel that passes through \(\tb{P}\)
and simplifying the action \(S\tb{P}S^{-1}\) of the spinor \(S=e^{-\tfrac{1}{2}\beta\mb{\Xi}}\).
Clifford translations can also be represented by quaternion multiplication.
If a point \(\tb{P}=w\e_{123}+x\e_{320}+y\e_{130}+z\e_{210}\) is normalised and
\(\mb{\Lambda}=\lambda_1\e_{23}+\lambda_2\e_{31}+\lambda_3\e_{12}\),
then the Clifford translation of \(\tb{P}\) by \(\beta\) generated by the positive Clifford bivector \(\mb{\Xi}=\mb{\Xi}^+=(\I+1)\mb{\Lambda}\)
is also given by \(p\mapsto pq\),
where \(p=w+x\tb{i}+y\tb{j}+z\tb{k}\) and \(q=\cos\beta-(\lambda_1\tb{i}+\lambda_2\tb{j}+\lambda_3\tb{k})\sin\beta\)
are the relevant quaternions.
So, a Clifford translation generated by \(\mb{\Xi}^+\) can also be referred to as a right Clifford translation.
For the negative Clifford bivector \(\mb{\Xi}=\mb{\Xi}^-=(\I-1)\mb{\Lambda}\),
 \(p\mapsto qp\) yields the same Clifford translation as that generated by \(\mb{\Xi}^-\).
It can be called a left Clifford translation.

\bibliographystyle{plain}
\bibliography{g.bib}

\end{document}